\documentclass[12pt,letterpaper,reqno]{amsart}

\usepackage{graphicx}

\usepackage[latin1]{inputenc}
\usepackage{amsmath, amsthm,  amsfonts, amscd}
\usepackage{amssymb} 
\usepackage{url}
\usepackage{hyperref}
\newtheorem{lemma}{Lemma}
\newtheorem{theorem}[lemma]{Theorem}
\newtheorem{corollary}[lemma]{Corollary}
\newtheorem{identity}[lemma]{Identity}
\newtheorem{conjecture}[lemma]{Conjecture}
\numberwithin{equation}{section} \numberwithin{lemma}{section}

\usepackage{cite}

\newcommand{\bigbox}[1]{
\begin{center}
\fbox{%
\begin{minipage}{\textwidth}
#1
\end{minipage}
}
\end{center}
}

\newcommand{\bbigbox}[1]{\bigbox{\textbf{#1}}}

\newcommand{\kommentar}[1]{}

\newcommand{\modd}{\text{~mod~}}
\newcommand{\RRe}{{\rm Re}}

\newcommand{\Neff}{N_{\rm eff}}

\newcommand{\FEp}{{\mathcal{F}}_{E}^+}

\newcommand*{\reff}[1]{(\ref{#1})}

\newcommand{\Nstd}{N_{\rm std}}

\DeclareMathOperator*{\Res}{Res}

\newtheorem{remark}[lemma]{Remark}
\newtheorem{proposition}[lemma]{Proposition}
\kommentar
{
\newtheorem{lemma}{Lemma}
\newtheorem{theorem}[lemma]{Theorem}

\numberwithin{equation}{section} \numberwithin{lemma}{section}
}

\kommentar
{

}

\def\thebibliography#1{\section{\centerline{\sc References}}
  \global\def\@listi{\leftmargin\leftmargini
               \labelwidth\leftmargini \advance\labelwidth-\labelsep
               \topsep 1pt plus 2pt minus 1pt
               \parsep 0.25ex plus 1pt \itemsep 0.25ex plus 1pt}
  \list {[\arabic{enumi}]}{\settowidth\labelwidth{[#1]}\leftmargin\labelwidth
    \advance\leftmargin\labelsep\usecounter{enumi}}
    \def\newblock{\hskip .11em plus .33em minus -.07em}
    \sloppy
    \sfcode`\.=1000\relax}

\begin{document}

\title[Random Matrix Model for Elliptic Curves with Finite Conductor]{A Random Matrix Model for Elliptic Curve $L$-Functions of Finite Conductor}

\vspace {2 in}

\author[Due\~nez]{E. Due\~nez}\email{eduenez@math.utsa.edu}
\address{Department of Mathematics, University of Texas at San Antonio, San Antonio, TX 78249, USA}

\author[Huynh]{D.K. Huynh}\email{dkhuynhms@gmail.com}
\address{Department of Pure Mathematics, University of Waterloo, Waterloo, ON, N2L 3G1, Canada}

\author[Keating]{J.P. Keating}\email{J.P.Keating@bristol.ac.uk}
\address{School of Mathematics, University of Bristol, Bristol BS8 1TW, UK}

\author[Miller]{S.J. Miller}\email{Steven.J.Miller@williams.edu}
\address{Department of Mathematics and Statistics, Williams College, Williamstown, MA 01267, USA}

\author[Snaith]{N.C. Snaith}\email{N.C.Snaith@bristol.ac.uk}
\address{School of Mathematics, University of Bristol, Bristol BS8 1TW, UK}



\date{\today}

\thanks{DKH was partially supported by EPSRC, a CRM postdoctoral fellowship and
  NSF grant DMS-0757627. SJM was partially supported by NSF grants DMS-0855257
  and DMS-097006. NCS was partially supported by funding from EPSRC. The work of
  JPK and NCS was sponsored by the Air Force Office of Scientific Research, Air
  Force Material Command, USAF, under grant number FA8655-10-1-3088. The
  U.S. Government is authorized to reproduce and distribute reprints for
  Government purpose notwithstanding any copyright notation thereon. DKH, JPK
  and NCS are grateful to the Mathematical Sciences Research Institute (MSRI),
  Berkeley, for hospitality when this work and manuscript were completed.
}

\begin{abstract}
  We propose a random matrix model for families of elliptic curve $L$-functions
  of finite conductor. A repulsion of the critical zeros of these $L$-functions away from the center
  of the critical strip was observed numerically by S. J. Miller in
  2006~\cite{Mil06}; such behaviour deviates qualitatively from the conjectural
  limiting distribution of the zeros (for large conductors this
  distribution is expected to approach the one-level density of eigenvalues of
  orthogonal matrices after appropriate rescaling).  Our purpose here is to
  provide a random matrix model for Miller's surprising discovery. We
  consider the family of even quadratic twists of a given elliptic curve. The
  main ingredient in our model is a calculation of the eigenvalue distribution
  of random orthogonal matrices
  whose characteristic polynomials are larger than some given value at the symmetry
  point in the spectra.
  We call this sub-ensemble of SO$(2N)$ the
  \emph{excised orthogonal ensemble.} The sieving-off of matrices with small
  values of the characteristic polynomial is akin to the discretization of the
  central values of $L$-functions implied by the formul\ae\ of Waldspurger and
  Kohnen-Zagier.  The cut-off scale appropriate to modeling elliptic curve
  $L$-functions is exponentially small relative to the matrix size~$N$. The
  one-level density of the excised ensemble can be expressed in terms of that of
  the well-known Jacobi ensemble, enabling the former to be explicitly
  calculated. It exhibits an exponentially small (on the
  scale of the mean spacing) hard gap determined by the cut-off value, followed
  by soft repulsion on a much larger scale. Neither of these features
  is present in the one-level density of $SO(2N)$. When $N \rightarrow \infty$ we recover the limiting orthogonal
  behaviour. Our results agree qualitatively with Miller's discrepancy. Choosing
  the cut-off appropriately gives a model in good quantitative agreement with
  the number-theoretical data.

\end{abstract}

\maketitle 


\section{Introduction} \label{introduction}

Understanding the ranks of elliptic curves is a central problem in
number theory. The celebrated Birch and Swinnerton-Dyer conjecture
provides an analytic approach to studying them via the order
of vanishing of the associated $L$-functions at the centre of the critical
strip: it states that the order of vanishing of the elliptic
curve $L$-function at the central point equals the rank of the
Mordell-Weil group. On the other hand, following the Katz-Sarnak philosophy, zero statistics of
families of $L$-functions are conjectured to have an underlying symmetry type associated with
certain random matrix ensembles~\cite{KatzSarnak99b,
KatzSarnak99a}: in an appropriate asymptotic limit, the statistics of the zeros of
families of $L$-functions are modeled by (a subgroup of) one of the classical
compact groups U$(N)$, O$(N)$ or USp$(2N)$.

If we could model the order of vanishing for families of elliptic curve
$L$-functions in terms of random matrices, then we would also gain information
about the distribution of ranks within a family of elliptic curves via the Birch
and Swinnerton-Dyer conjecture.  This is the wider goal that motivates the
current work.

There has been much success modelling families of $L$-functions with matrices
from the classical compact groups
\cite{BuiKeat1,BuiKeat2,CRS05,ConSna07,ConSna08,ConSna08a,ConFarm00,CFKRS05,CFZa,CKRS02,CKRS06,CRSW05,duenez2006low,FGH05,FouIwa03,GHK07,GoMil,Gul05,Hug05,HKO00,Hug03,
  HugMil07,HugRud03,HugRud02,HugRud03b,ILS99,KatzSarnak99b,KatzSarnak99a,KeatSnaith00a,KeatSnaith00b,
  Mez03,MilPe,kn:ozlsny99, RicottaRoyer09,royer2001petits,Rub01,Sna05},
including families of elliptic curve
$L$-functions~\cite{CKRS02,CKRS06,HKS,HMM,MillerPhD,Mil04_a,Young06}.
This work strongly supports the Katz-Sarnak philosophy that, in
the correct asymptotic limit, matrices from the classical compact
groups accurately model the zero statistics of $L$-functions.
Specifically, for families of elliptic curve $L$-functions it has
been shown for restricted test-functions that the one- and
two-level densities do indeed show the expected symmetry of
orthogonal type\footnote{\label{foot:ellcurverank} If the family
has geometric rank $r$,
  then by the Birch and Swinnerton-Dyer conjecture and Silverman's
  specialization theorem all but finitely many of the elliptic curve
  $L$-functions have $r$ zeros at the central point; the correct scaling limit
  is the limit of block diagonal orthogonal matrices with an $r\times r$
  identity matrix as the upper left block.}, that is, of random matrices from
O$(N)$, in the limit of large conductor~\cite{MillerPhD, Mil04_a, Young06}.

For the Riemann zeta function, and more generally families of $L$-functions,
good qualitative agreement is also found, as the asymptotic limit is approached,
by modeling the $L$-functions with finite-sized matrices where the matrix size
is chosen so that the mean density of matrix eigenvalues matches the mean
density of zeros of the $L$-functions under consideration (originally proposed
by Keating and Snaith~\cite{KeatSnaith00a,KeatSnaith00b}). This finite matrix
size model has been refined still further with the incorporation of
number-theoretic information, leading to extremely convincing models for, for
example, correlation functions~\cite{BerryKeat99b, BogKeat96b}, spacing
distributions~\cite{BBLM06}, and moments and ratios of values of
$L$-functions~\cite{CFKRS05,CFZa}.

However, in the elliptic curve case S. J. Miller, in numerical
data gathered in~\cite{Mil06}, discovered a significant
discrepancy from the scaling limit of the expected model of
orthogonal matrices.  That is, for elliptic curve $L$-functions of
finite conductor there is a large (repulsive) deviation of the
zero statistics from the expected orthogonal symmetry that is not
explained by either using finite-size matrices or by the inclusion
of arithmetical terms as in \cite{BerryKeat99b, BBLM06,
BogKeat96b, CFKRS05, CFZa}.  Miller was the first to discover the
soft repulsion of zeros of elliptic curve $L$-functions from the
central point, and his data is reproduced in
figure~\ref{fig:miller}.  However, the larger data sets of the
present paper (see figure~\ref{fig:p6_vs_firstzeros}), along with
the random matrix model, indicate a hard gap containing no
zeros.\footnote{Subsequent to our posting this manuscript on the
arXiv, Marshall \cite{kn:marshall} has developed a rigorous theory
for this hard gap.} The specific goal of this paper is to present
a model that accurately describes the zero statistics of elliptic
curve $L$-functions of finite conductor.

\newcommand{\XX}{{\mathcal X}}
\newcommand{\CSO}{C_{{\rm SO}(2N)}}
\newcommand{\CX}{C_{\XX}}
\newcommand{\RTX}{R_1^{T_\XX}}
\newcommand{\TX}{T_{\XX}}
\newcommand{\SO}{{\rm{SO(2$N$)}}\/}
The main ingredient in the new model is to mimic central values of elliptic curve $L$-functions by the value of the characteristic polynomial at 1 for matrices from SO$(2N)$. The key fact about these central values, not previously incorporated into random matrix models for the zero statistics, is that
they are either zero (for curves with rank greater than zero) or they are greater or equal to a finite minimal size fixed by arithmetical considerations.
The characteristic polynomial of a matrix $A\in$ SO$(2N)$ is given by
\begin{equation}
\Lambda_A(e^{i\theta},N) := \det(I-e^{i\theta}A^{-1}) =
\prod_{k=1}^N (1 - e^{i(\theta - \theta_k)})(1 - e^{i(\theta +
\theta_k)}),
\end{equation}
with $e^{\pm i\theta_1},\ldots,e^{\pm i\theta_N}$ the eigenvalues
of $A$. Motivated by the arithmetical size constraint on the
central values of the $L$-functions we seek to understand, we here
consider the set $\TX$ of matrices $A \in$ SO$(2N)$ whose
characteristic polynomial $\Lambda_A(1,N)$ is larger than some
cut-off value $\exp(\XX)$. We call this the \emph{excised
orthogonal ensemble}.  In figure~\ref{fig:SO2Nvsexcised} we plot
the distribution of the first eigenvalue of $SO(24)$, which shows
no repulsion, as well as the distribution of the first eigenvalue
of an excised ensemble of $24\times 24$ orthogonal matrices.  Here
the size of the excision is $\exp(\XX)\approx 0.005$ and we see
the hard gap at the origin where no eigenvalues lie, as well as
the soft repulsion that interpolates between zero and the bulk of
the distribution.  This is seen more clearly in the right hand
figure where a region near the origin is enlarged.

\begin{figure}[h]
\begin{center}
 \includegraphics[width=6.48cm]{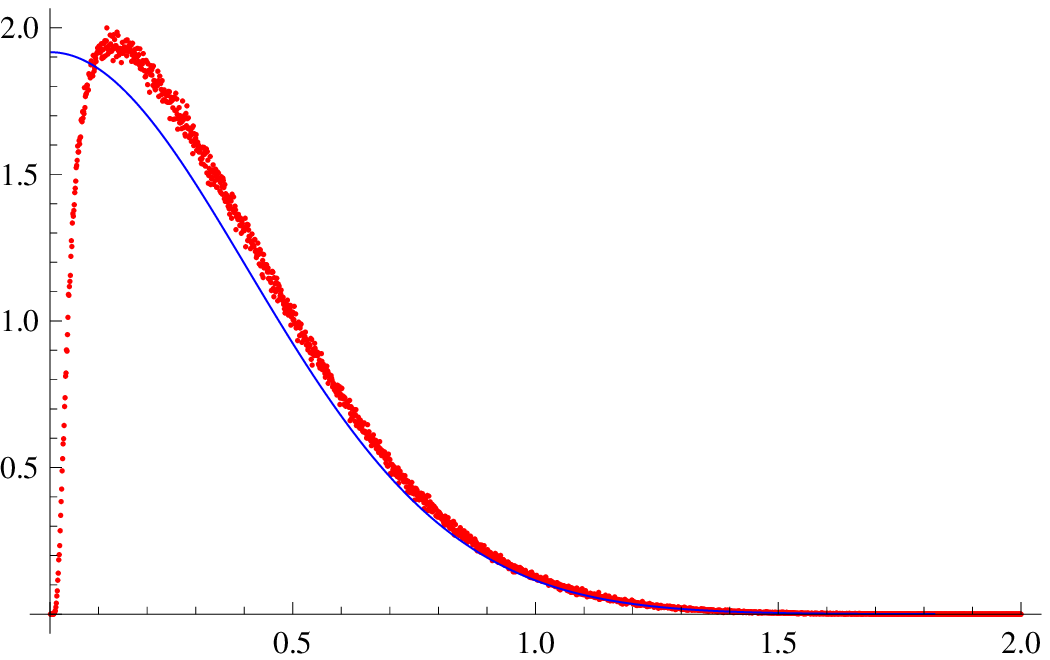}\ \includegraphics[width=6.48cm]{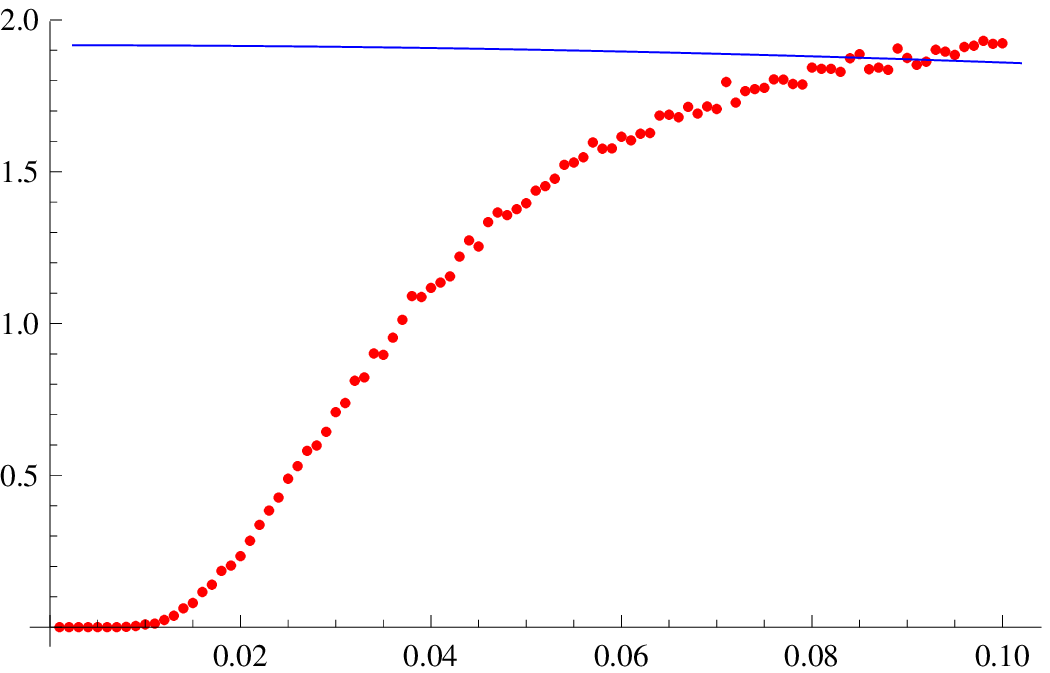}
\caption{\label{fig:SO2Nvsexcised} The distribution (blue, solid
curve) of the first eigenvalue of $SO(24)$, showing no repulsion,
and the distribution (red dots) of the first eigenvalue of an
excised ensemble of $24\times 24$ orthogonal matrices with
$\exp(\XX)\approx 0.005$. The right hand figure shows an
enlargement near the origin. The $SO(24)$ curve is computed using
the numerical differential equation solver for Painlev\'e~VI which
is developed in~\cite{DHKMS2}.  For the excised ensemble, 3 000
000 random $SO(24)$ matrices were generated and those not
satisfying $\Lambda_A(1,N)\geq\exp(\XX)$ were discarded.}
\end{center}
\end{figure}

Recall that the one-level density $R^{G(N)}_1$ for a (circular) ensemble $G(N)$
of matrices whose eigenvalues are parametrized by an unordered $N$-tuple of
eigenphases $\{\theta_n\}_1^N$ is given by
\begin{equation}\label{eq:genonelevel}
R^{G(N)}_1(\theta) = N \int \ldots \int P(\theta,\theta_2,\ldots,\theta_N) d\theta_2\ldots d\theta_N
\end{equation}
where $P(\theta,\theta_2,\ldots,\theta_N)$ is the joint probability density
function of eigenphases\footnote{Note that in the case of our interest, namely
  the (full, or excised) ensemble of orthogonal matrices of size~$2N$, there are
  $N$ pairs of eigenvalues parametrized by $N$ eigenphases.}.  Helpfully, the
probability density function of the one-level density $\RTX$ can be expressed in
terms of the well-known Jacobi ensemble $J_N$ (see \cite{LGRM} for properties of this ensemble):
\begin{theorem} \label{lemmaone}
The one-level density $\RTX$ for the set $\TX$ of matrices $A \in$ SO$(2N)$
with characteristic polynomial $\Lambda_A(e^{i\theta},N)$ satisfying $\log|\Lambda_A(1,N)| \geq \XX$ is
given by
\begin{align*}
R_1^{\TX}(\theta_1) & = \frac{\CX}{2\pi i} \int_{c-i\infty}^{c+i\infty}2^{Nr}\frac{\exp(-r\XX)}{r} R_1^{J_N}(\theta_1; r-1/2,-1/2)dr
\end{align*}
where $\CX$ is a normalization constant defined in~\reff{condition_on_CX} and
\begin{multline}\label{eqn:integrand}
  R_1^{J_N}(\theta_1; r-1/2,-1/2)  = N \int_0^\pi\cdots\int^\pi_0 \prod_{j=1}^N w^{(r-1/2,-1/2)}(\cos \theta_j) \\
  \times \prod_{j<k}(\cos \theta_j-\cos \theta_k)^2d\theta_2 \cdots d\theta_N
\end{multline}
is the one-level density for the Jacobi ensemble $J_N$ with weight function
\begin{equation*}
w^{(\alpha, \beta)}(\cos \theta) = (1-\cos \theta)^{\alpha + 1/2} (1+\cos \theta)^{\beta + 1/2},
\qquad\text{$\alpha = r - 1/2$ and $\beta = -1/2.$}
\end{equation*}
\end{theorem}

Applying the method of orthogonal polynomials we evaluate $R_1^{J_N}(\theta;\alpha,\beta)$ to give
\begin{theorem} \label{theoremone}
With $\RTX$ as above, we have
\begin{equation} \label{eq:theoremone}
  \begin{split}
    R_1^{\TX}(\theta) &= \frac{\CX}{2\pi i} \int_{c-i\infty}^{c+i\infty}
    \frac{\exp(-r\XX)}{r} 2^{N^2+2Nr-N} \times\\
    &\times \prod_{j=0}^{N-1} \frac{\Gamma(2+j)\Gamma(1/2+j)\Gamma(r+1/2+j)}
    {\Gamma(r+N+j)} \times\\
    &\times (1-\cos \theta)^r \frac{2^{1-r}}{2N+r-1}
    \frac{\Gamma(N+1)\Gamma(N+r)} {\Gamma(N+r-1/2)\Gamma(N-1/2)}
    P(N,r,\theta)\,dr
  \end{split}
\end{equation}
with normalization constant $\CX$ defined in~\reff{condition_on_CX} and $P(N,r,\theta)$ defined
in terms of Jacobi polynomials in~\reff{PNrtheta}.
\end{theorem}

An immediate consequence using residue calculus is the following:
\begin{theorem} \label{corollaryone}
With $\RTX$ as above, we have
\begin{equation} \label{eq:corollary}
\RTX(\theta) =
\begin{cases}
0 & \mbox {~for~} d(\theta, \XX) < 0\\
R_1^{{\rm SO}(2N)}(\theta) + \CX \sum_{k=0}^\infty b_k(\theta)
\exp((k+1/2) \XX) & \mbox{~for~} d(\theta, \XX) \geq 0,
\end{cases}
\end{equation}
where $d(\theta, \XX):= (2N-1)\log2 +\log(1-\cos\theta)-\XX$ and
$b_k$ are coefficients arising from the residues. The
normalization constant $\CX$ is defined in~\reff{condition_on_CX}
and $R_1^{{\rm SO}(2N)}(\theta)$ is the one-level density for
$SO(2N)$, defined by  (\ref{eq:genonelevel}).
\end{theorem}
\noindent Thus in the limit $\XX \rightarrow -\infty$, $\theta$
fixed, $R_1^{\TX}(\theta) \rightarrow R_1^{{\rm SO}(2N)}(\theta)$.

In order to compare our random matrix results with the elliptic
curve data we shall need to take $-\XX$ proportional to $N$, i.e.
the cut-off is exponentially small in $N$.  From
Theorem~\ref{corollaryone}, and illustrated in
figure~\ref{fig:consistency}, one sees the existence of a hard gap
where $\RTX(\theta)$ vanishes, namely the interval $\{\theta>0
\mid d(\theta,\XX)<0\}$.  When $-\XX \propto N$ the hard gap is
therefore exponentially small in $N$. Beyond this range the
formula exhibits soft repulsion of the eigenvalues from 1 on a
much larger scale; the repulsion extends  far beyond the hard gap.
It is interesting that such a `pierced' subset of SO$(2N)$, where
the cut-off is exponentially small on the scale of the mean
eigenvalue spacing, has such a pronounced effect on the eigenvalue
statistics over a much larger distance.

This agrees qualitatively with the discrepancy Miller observed.
We go on to discuss two ways to determine the relevant parameters which lead not just to qualitative but also quantitative agreement.

In the following sections we give some background information on elliptic curves
(section~\ref{sect:background}), detail the anomalous zero
statistics observed by Miller in~\cite{Mil06} (section
\ref{sect:miller}) and set out the new model (section
\ref{sect:model}). The model involves selecting from SO$(2N)$ those matrices whose characteristic polynomial at 1 is larger than some given cut-off value (section~\ref{sect:discvalue}).  We apply the model with two matrix sizes: one the standard value related to the mean density of zeros and one an effective matrix size
 determined from lower-order terms of the
one-level density (section~\ref{sect:neff}). We present numerical
evidence for our model using the distribution of the first zero of
a family of elliptic curve $L$-functions (sections
\ref{sect:discevidenceN0},~\ref{sect:discevidenceNeff}) and also
using the one-level density statistic for the excised ensemble of
matrices (section~\ref{sect:onelevel}).


\section{Elliptic curve $L$-functions}

\subsection{Background}\label{sect:background}

An elliptic curve $E$ can be written in Weierstrass form as
\begin{equation}
\label{eq:ec}  y^2+c_1xy+c_3y=x^3+c_2x^2+c_4x+c_6, \;\; c_i\in
\mathbb{Z}.
\end{equation}
The $L$-function $L_E(s)$ associated with $E$ is given by the
Dirichlet series
\begin{equation}
L_E(s)=\sum_{n=1}^{\infty} \frac{\lambda(n)}{n^s},
\end{equation}
where the coefficients ($\lambda(n)=a(n)/\sqrt{n}$, with
$a(p)=p+1-\#E({\mathbb F}_p)$, $\#E({\mathbb F}_p)$ being the
number of points on $E$ counted over ${\mathbb F}_p$), have been
normalized so that the functional equation relates $s$ to
$1-s$:
\begin{equation}
 L_E(s)=\omega(E) \left(\frac{2\pi}{\sqrt{M}}\right)^{2s-1}
\frac{\Gamma(3/2-s)}{ \Gamma(s+1/2)}L_E(1-s).
\end{equation}
Here $M$ is the conductor of the elliptic curve $E$; for convenience we will
consider only prime $M$. Also, $\omega(E)$ is $+1$ or $-1$
resulting, respectively, in an even or odd functional equation for
$L_E$.  There is a Generalized Riemann Hypothesis for
$L$-functions of  elliptic curves, stating that the non-trivial
zeros of $L_E(s)$ lie on the critical line (with real part equal
to 1/2).  The even (odd) symmetry of the functional equation
implies an even (odd) symmetry of the zeros around the central
point $s = 1/2$, where the critical line crosses the real axis.  Often when
we consider zeros statistics of a family of $L$-functions we are
particularly interested in the zeros close to the central point.

The family of elliptic curve $L$-functions for which numerical
evidence is presented in this paper is that of quadratic twists of
a fixed curve $E$.  Let $L_E(s, \chi_d)$ denote the $L$-function
obtained by twisting $L_E(s)$ by a quadratic character $\chi_d$.
Here $d$ is a fundamental discriminant, i.e., $d \in
\mathbb{Z}$ such that $p^2 \nmid d$ for all odd primes $p$ and
$d \equiv 1 \modd 4$ or $d \equiv 8, 12 \modd 16$, and $\chi_d$ is
the Kronecker symbol (an extension of the Legendre symbol, taking
the values 1, 0, or $-1$). The twisted $L$-function, which is itself
the $L$-function associated with another elliptic curve $E_d$, is
given by
\begin{equation} \label{twistedLfunction}
L_E(s, \chi_d) = \sum_{n = 1}^\infty \frac{\lambda(n)
\chi_d(n)}{n^s} = \prod_p \left(1 -
\frac{\lambda(p)\chi_d(p)}{p^s} +
\frac{\psi_{M}(p)\chi_d(p)^2}{p^{2s}} \right)^{-1},
\end{equation}
where $\psi_{M}$ is the principal Dirichlet character of modulus
$M$:
\begin{equation}\label{psi}
\psi_{M}(p) =
\begin{cases}
1 \text{~if~} p \nmid M\\
0 \text{~otherwise}.
\end{cases}
\end{equation}
The functional equation of this $L$-function is, for $(d,M)=1$,
\begin{equation}
L_E(s,\chi_d)=\chi_d(-M)\omega(E)
\left(\frac{2\pi}{\sqrt{M}|d|}\right)^{2s-1} \frac{\Gamma(3/2-s)}{
\Gamma(s+1/2)}L_E(1-s,\chi_d).
\end{equation}

The sign of this functional equation is $\chi_d(-M)\omega(E)$ and
it is more instructive to restrict to the fundamental
discriminants which, for a fixed $E$, give an even, or
alternatively an odd, functional equation for $L_E(s,\chi_d)$.  In
particular, the families we will consider in this paper will be
denoted $\FEp(X)$: those quadratic twists of the curve $E$  that
have an even functional equation and $0<d\leq X$ (or else $-X\leq
d<0$). In the appropriate asymptotic limit, zero statistics of
families of $L$-functions of elliptic curves with even functional
equation are believed to follow the distribution laws of
eigenvalues of the even orthogonal group SO$(2N)$ and those with
odd functional equation are expected to show SO$(2N+1)$
statistics.\footnote{If the family has rank then we must modify the random matrix ensemble as in Footnote~\ref{foot:ellcurverank}.}  The asymptotic parameter of the family $\FEp(X)$ is
$X$.  The conductor of $L_E(s,\chi_d)$ is $Md^2$, if $M$ is the
conductor of the original curve, and so $d$ is the parameter that
orders the curves by conductor.  It is expected that as
$X\rightarrow \infty$ the zero statistics of $\FEp(X)$ tend to the
large-$N$ limiting statistics of eigenvalues from SO$(2N)$ in
accordance with the Katz-Sarnak philosophy. We propose a model
for the behaviour of zero statistics for finite $X$.

Key to this model are formul\ae\ of Waldspurger~\cite{Walds80},
Kohnen-Zagier~\cite{KohnenZagier81} and Baruch-Mao
\cite{BaruchMao2007} which show that modular $L$-functions can
only attain discrete values at the centre of the critical
strip. In particular we have for the twists of an elliptic curve
$L$-function
\begin{equation} \label{Waldspurgerformula}
L_E(1/2, \chi_d) = \kappa_E \frac{c_E(|d|)^2}{d^{1/2}},
\end{equation}
where $c_E(|d|)$ are integers and the Fourier coefficients of a weight-$3/2$
modular form and $\kappa_E$ is a constant depending on the curve $E$.

\subsection{Unexpected repulsion}\label{sect:miller}

\begin{figure}[h]
\begin{center}
 \includegraphics[width=6.48cm]{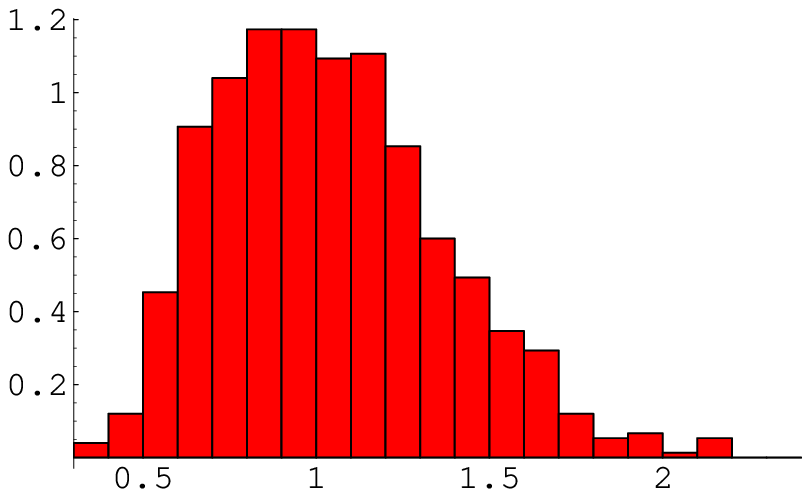}\ \includegraphics[width=6.48cm]{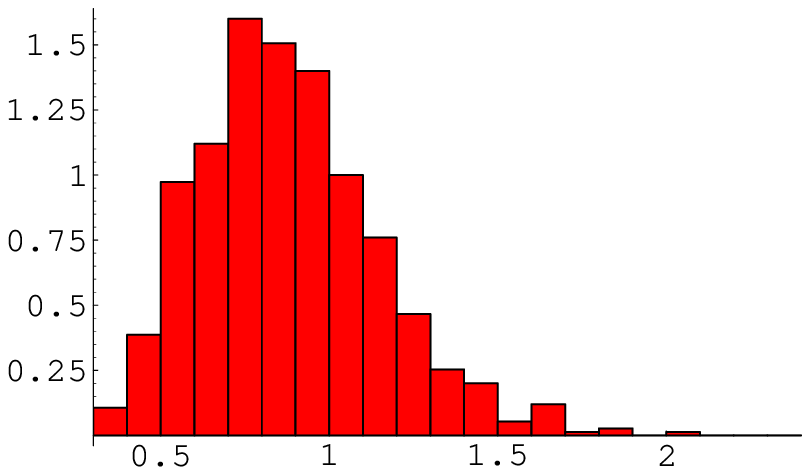}
\caption{\label{fig:miller} First normalized zero above the
central point: Left: $750$ rank 0 curves from $y^2 +
a_1xy+a_3y=x^3+a_2x^2+a_4x+a_6$, $\log({\rm cond}) \in [3.2,
12.6]$, $\text{median} = 1.00$, $\text{mean} = 1.04$, standard
deviation about the mean $=.32$. Right: $750$ rank 0 curves from
$y^2 + a_1xy+a_3y=x^3+a_2x^2+a_4x+a_6$, $\log({\rm cond}) \in
[12.6, 14.9]$, $\text{median}=.85$, $\text{mean} = .88$, standard
deviation about the mean $= .27$}
\end{center}
\end{figure}

In 2006 S. J. Miller~\cite{Mil06} investigated the statistics of
the zeros of various families of elliptic curve $L$-functions.
Figure~\ref{fig:miller} gives an example of what he discovered. It
shows a histogram of the first zero above the central point for
rank zero elliptic curve $L$-functions generated by randomly
selecting the coefficients $c_1$ up to $c_6$ (as defined in
\eqref{eq:ec}) for curves with conductors in the ranges indicated
in the caption. The zeros are scaled by the mean density of low
zeros of the $L$-functions and the plots are normalized so that
they represent the probability density function for the first zero
of $L$-functions from this family. Miller observes that there is
clear repulsion of the first zero from the central point; that is,
the plots drop to zero at the origin, indicating a very low
probability of finding an $L$-function with a low first zero.

\begin{figure}[h]
\begin{center}
 \includegraphics[width=10cm]{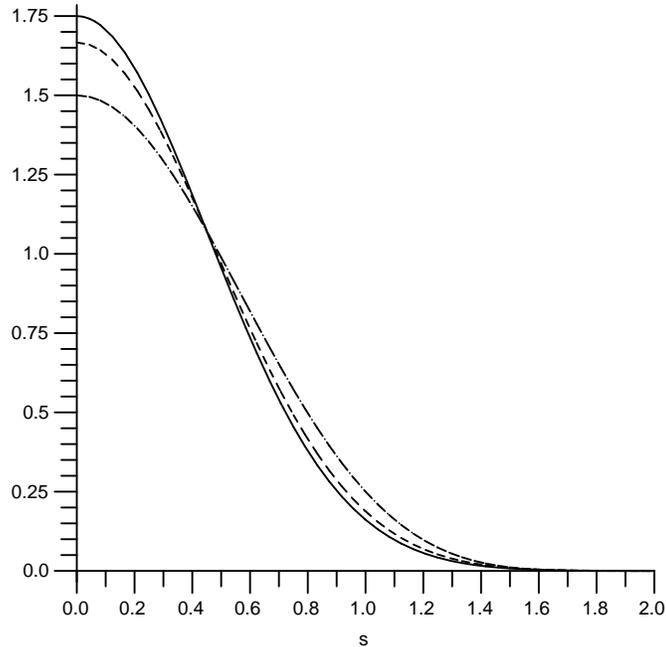}
\caption{\label{fig:millerRMT} Probability density of normalized
eigenvalue closest to 1 for SO$(8)$ (solid), SO$(6)$ (dashed)
and SO$(4)$ (dot-dashed).}
\end{center}
\end{figure}

What is surprising about these plots is that the standard way to model such
$L$-functions would be with matrices from SO$(2N)$, with $N$ chosen to be equal
to half the logarithm of the conductor of the curve.  This choice of $N$ has the
effect of equating the density of eigenvalues to the density of zeros near the
critical point and there has been much work showing, for the Riemann zeta
function~\cite{KeatSnaith00a,GHK07,FGH05,HKO00,Hug03,
  HugRud02,Mez03,CRS05,ConSna07,ConSna08a,ConSna08} and for families of
$L$-functions not necessarily associated with elliptic curves
\cite{KeatSnaith00b,ConFarm00,CFKRS05,CFZa,CKRS02,CKRS06,CRSW05,BuiKeat1,
  BuiKeat2,HugRud03b,Sna05,ConSna07}, that this can be effective at modelling
number-theoretic data even far from the asymptotic limit specified by Katz and
Sarnak. However, as figure~\ref{fig:millerRMT} illustrates, even for small-size
SO$(2N)$ matrices, the distribution of the eigenvalue closest to 1 on the unit
circle (the random matrix equivalent to the distribution of the first zero of
$L$-functions) shows no repulsion at the origin of the distribution ---a fact of
course well-known in random matrix theory.

This, then, is the mystery.  No one doubts that in the large-conductor limit the
distribution of the first zero of elliptic curve $L$-functions will tend to the
large-$N$ limit of the distribution of the first eigenvalues of SO$(2N)$
matrices: Miller observes (see figure~\ref{fig:miller}) that the repulsion
decreases with increasing conductor.  However, for finite conductor we do not
get qualitative agreement with the statistics of finite-sized matrices, as we do
in many other families of $L$-functions. Indeed one sees a qualitatively
unexpected phenomenon, namely, repulsion of zeros away from the central point.
In section~\ref{sect:model} we outline an idea that successfully models the
anomalous statistical behaviour discovered by Miller.

In this paper we concentrate on developing a model for the
unexpected zero statistics of rank zero curves seen in figure
\ref{fig:miller}.   However, it should be noted that Miller
\cite{Mil06} also investigated zeros of higher rank curves and the
behaviour of zeros further from the central point, which will
guide future attempts to model statistics of $L$-functions
associated to elliptic curves of higher rank. Briefly, the
numerical findings were the following:

\begin{enumerate}

\item The repulsion of the low-lying zeros increased with increasing rank, and
  was present even for rank 0 curves.

\item As the conductors increased, the repulsion decreased.

\item Statistical tests failed to reject the hypothesis that, on average, the
  first three zeros were all repelled equally (i.~e., shifted by the same
  amount).
\end{enumerate}


\section{The model}\label{sect:model}

We propose to model the zero statistics of rank-zero elliptic
curve $L$-functions with the subset $\TX$ of $SO(2N)$ defined in
the introduction. The justification for employing an
\emph{excised} random matrix model such as this is that formula
(\ref{Waldspurgerformula}) indicates that the $L$-functions
themselves have a discretization at the central point.  Thus the
statement
\begin{equation} \label{discretizationBogo}
L_E(1/2, \chi_d) < \frac{\kappa_E}{d^{1/2}}
\end{equation}
implies that
\begin{equation}
L_E(1/2, \chi_d) = 0.
\end{equation}
Hence, each $L_E(s, \chi_d) \in \FEp$ associated with a rank 0
curve $E_d$ satisfies
\begin{equation} \label{rankzerodiscr}
L_E(1/2, \chi_d) \geq \frac{\kappa_E}{d^{1/2}}.
\end{equation}
Thus, since previous
work~\cite{KeatSnaith00b,ConFarm00,CKRS02,CKRS06} shows that
values of $L$-functions at the central point can be characterized
using characteristic polynomials evaluated at the point 1, to
model rank 0 curves we discard from our orthogonal ensemble all
matrices not satisfying (\ref{eq:genericdisc}) and we will propose
a value of $\XX$ based on (\ref{rankzerodiscr}).  For a start, if
we are working with discriminants of size around $d$ then equating
the density of eigenvalues, $N/\pi$, with that of zeros near the
central point, $\frac{1}{\pi} \log (\tfrac{\sqrt{M} d}{2\pi})$,
gives us an equivalent value of $N$: $\Nstd\sim\log d$.  Since the
$L$-values are discretized on a scale of $1/\sqrt{d}$, we
``excise" (i.~e., discard), characteristic polynomials whose value
at $1$ is of the scale~$\exp(-\Nstd/2)$.

Thus, for certain values of $N$ we plan to model elliptic curve zero statistics using matrices from SO$(2N)$ satisfying
\begin{equation}\label{eq:genericdisc}
|\Lambda_A(1,N)|\geq \exp{\XX}=c \times \exp(-\Nstd/2),
\end{equation}
and we will propose a value for $c$ in section~\ref{sect:discvalue}.  We present the model for the zero statistics using two different matrix sizes. One obvious choice is $N=\Nstd$ and we present data for this case in section~\ref{sect:discevidenceN0}.  The other choice is an ``effective" matrix size, $\Neff$, that incorporates arithmetic information and is based on an idea of Bogomolny, Bohigas, Leboeuf
and Monastra~\cite{BBLM06}.  We illustrate how $\Neff$ is calculated in section~\ref{sect:neff} and present the numerical results of the model in section~\ref{sect:discevidenceNeff}.  In this paper we compare this model with numerical data from a family of twists of one particular elliptic curve, $E_{11}$.


\section{lower-order terms and the effective matrix
size}\label{sect:neff}

By the Katz-Sarnak philosophy the statistical properties of zeros of
$L$-functions should asymptotically behave like the scaling limits of
eigenvalues of random matrices drawn from one of the classical compact
groups. However, for finite values of the asymptotic parameter we observe
deviations (of number-theoretic origin) from the limiting result.  We note that
for many families of $L$-functions these deviations, at least over length scales
on the order of the mean spacing of zeros, are slight perturbations of the
limiting results, in contrast to the elliptic curve case where we see distinct,
qualitative disagreement (the repulsion at the central point).\footnote{Random
  matrix theory does not see the fine arithmetic properties of the family, which
  surface in the lower-order terms.  Thus, while the main terms of various
  families of elliptic curves are the same, the lower-order terms show
  differences due, for instance, to complex multiplication or torsion points;
  see~\cite{Mil09a}. }

In this section we adapt the method of Bogomolny, Bohigas, Leboeuf
and Monastra~\cite{BBLM06}, previously used to investigate statistics of
zeros of the Riemann zeta function at a finite height on the
critical line, and apply it to the family $\FEp(X)$ of elliptic
curve $L$-functions.

Specifically, Bogomolny \emph{et al.} use a conjectured
formula for the two-point correlation function of the Riemann
zeros at finite height on the critical line~\cite{BogKeat96b} and
compare it to the finite-$N$ form of the two point correlation
function of unitary matrices from U$(N)$.  Under the standard
procedure of equating mean densities of eigenvalues with the mean
density of zeros, yielding in this case
\begin{equation} N =
\log\left(\frac{T}{2\pi}\right),
\end{equation}
the leading order term of these two formul\ae\ match up, but
Bogomolny \emph{et al.} show that by looking at a scaled version of the
statistic and then choosing an ``effective" matrix size, related
to $N$ by multiplication by a constant of arithmetic origin, they
can match the first lower-order term in the formul\ae.  By scaling
a further variable they match yet another lower-order term, although this step requires that the considered Riemann
zeros be high lying, a point we will return to later.  The authors illustrate with comprehensive numerical results that with
each successive refinement the fit of the random-matrix model to
the Riemann zero data significantly improves.  What is remarkable
about this is that they don't just see good numerical agreement in
the two-point correlation statistic, where good agreement is to be expected as terms were matched by design, but
also in the nearest-neighbour spacing distribution.  The nearest
neighbour spacing statistic is the probability density for
distances between consecutive zeros, or equivalently, a normalized
histogram of gaps between consecutive zeros.  The argument is
that the corrections to the asymptotics in the form of $N_{\rm
eff}$ and the further scaling factor is also valid for all
correlation functions, and therefore also for the nearest
neighbour spacing.  The strength of this work is that using this heuristic approach the information gained from a simpler statistic (two-point correlation) yields information
for a more complicated one (nearest-neighbour spacing)---the latter
one being determined by all correlation functions together.

\kommentar
{
\textbf{(COMMENT FROM SJMILLER: I WORKED ON THIS WITH A THESIS STUDENT, RALPH MORRISON, LAST YEAR AT WILLIAMS. WE PROVED IT FOR THE TAU FUNCTION INSTEAD OF AN ELLIPTIC CURVE -- THIS REMOVED THE LEVEL ASPECT. OUR ARGUMENTS SHOULD APPLY TO OUR CASE HERE AS WELL -- WE'VE ALMOST FINISHED THIS PROOF, HOPEFULLY JUST A FEW DAYS AWAY, SEE~\cite{HMM}.)}
}
In the following we adopt the above method and apply it to the
case of $\FEp(X)$.  In this case we have a conjectural formula for
the one-level density for the scaled zeros of $\FEp(X)$, as derived
in~\cite{HKS,HMM}.
For our purposes this is best given as the
following expansion for large $X$ (equation (3.18) in~\cite{HKS}):
\begin{multline}
\frac{1}{X^*}\sum_{\substack{0<d \leq X\\
\chi_d(-M)\omega_E=+1}} \sum_{\gamma_d}g\Big(\frac{\gamma_d L}{\pi}\Big) \\
 =\int_{-\infty}^{\infty}g(\tau)\Bigg( 1 + \frac{\sin(2\pi
\tau)}{2 \pi \tau} - r_1\frac{1 +\cos(2\pi \tau)}{L} -
r_2\frac{\pi \tau \sin(2\pi \tau)}{L^2}\Bigg)d\tau +
O\left(\frac{1}{L^3}\right)
\label{S1f}
\end{multline}
with
\begin{equation}
L = \log\bigg(\frac{\sqrt{M}X}{2\pi} \bigg),
\end{equation}
where $\gamma_d$ is the imaginary part of a generic zero of
$L_E(s, \chi_d)$, the sum is over fundamental discriminants $d$,
and $X^*$ is the number of fundamental discriminants satisfying
the conditions on the sum. The coefficients $r_1$ and $r_2$ in
\reff{S1f} are arithmetic constants involving the Dirichlet
coefficients of $L_E(s)$; see~\cite{HKS} for details.

As expected the result~\reff{S1f} has the form of the asymptotic
random-matrix result
\begin{equation} \label{R1SO2N}
\widetilde{R}_1(s) = 1 + \left(\frac{\sin 2\pi s}{2 \pi s}\right)
\end{equation}
for the even orthogonal group, plus correction terms.  We compare
this expansion with the one we obtain by expanding the scaled
one-level density of SO$(2N)$ for finite $N$. The unscaled one-level
density of SO$(2N)$ is  (see, for example,~\cite{Con05a})
\begin{equation} \label{R1SO2Nunscaled}
R_1(s) = \frac{2N-1}{2\pi} + \frac{\sin((2N-1)s)}{2\pi \sin s},
\end{equation}
and so scaling by the mean density and expanding in powers of $1/N$ gives
\begin{equation} \label{SO2N}
\frac{\pi}{N} R_1\left(\frac{\pi y}{N}\right) = 1 +
\frac{\sin(2\pi y)}{2\pi y} - \frac{1 + \cos(2\pi y)}{2N} -
\frac{\pi y \sin(2\pi y)}{6N^2} + O\left(\frac{1}{N^3}\right).
\end{equation}

By choosing an effective matrix size
\begin{equation}\label{eq:defNeff}
\Neff = \frac{L}{2r_1}
\end{equation}
we match the next-to-leading term in (\ref{S1f}) and (\ref{SO2N}). Arguing as
Bogomolny \emph{et al.} we conjecture that the improvement made by using
matrices of size $\Neff$ also holds for all $n$-point correlation or density
functions.  In particular we apply this result to the distribution of the lowest
zero and see significantly better agreement in the bulk and tail of the
distribution when we use $\Neff$ as opposed to the use of the standard matrix
size $\Nstd = L$.

\begin{figure}[h]
\includegraphics[scale=.50, angle=-90]{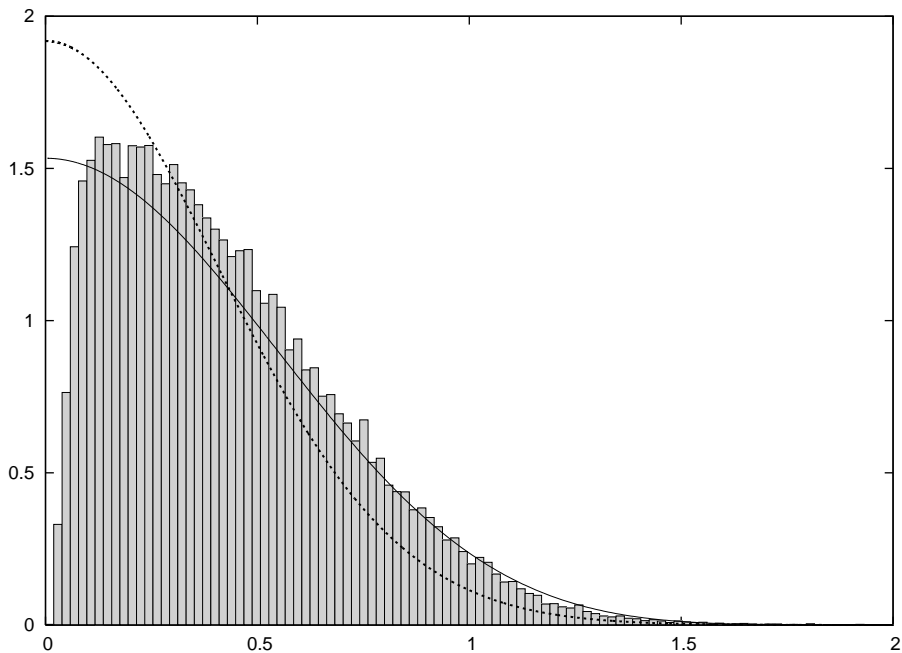}
\caption{Distribution of the lowest zero for $L_{E_{11}}(s,
\chi_d)$ with $0 < d \leq \text{400,000}$ (bar chart), distribution of
the lowest eigenvalue of SO$(2N)$ with $N_{\rm eff}=2.14$ (solid),
standard $\Nstd=12.26$ (dots). }\label{fig:p6_vs_firstzeros}
\end{figure}

We illustrate the effect of $\Neff$ in figure
\ref{fig:p6_vs_firstzeros}. We choose the elliptic curve $E_{11}$
($(c_1,c_2,c_3,c_4,c_6)=(0,-1,1,0,0)$ in Weierstrass form)
and find numerically that, for this curve,
\begin{align}\label{eq:r1E11}
r_1 & \approx 2.8600.
\end{align}
\kommentar
{
\bbigbox{**** Duc Khiem, do you have a numerical value for $r_2$?  If so,
let's state it here so that it could be used by anyone who could
compute more data than we could. *****}
}

We used Rubinstein's {\sffamily lcalc}~\cite{Rub} to compute the lowest zero of
the quadratic twists $L_{E_{11}}(s,\chi_d)$ with even functional equation and
with fundamental discriminants $d$, $0 < d \leq X = \text{400,000}$. The
standard matrix size corresponding to $X = \text{400,000}$ is
\begin{equation}
\Nstd = \log \left(\frac{\sqrt{11} X}{2\pi} \right) \approx 12.26,
\end{equation}
whereas the effective one is
\begin{equation}
\Neff = \frac{\Nstd}{2r_1} \approx 2.14.
\end{equation}
The distribution of the lowest zero (a normalized histogram of the heights of
the lowest zero of each $L$-function in the family) is then depicted (for
rank-zero curves) as a bar chart in figure \ref{fig:p6_vs_firstzeros} whereas
the distribution of the lowest eigenvalue of SO$(2N)$ with `effective' matrix
size $N_{\rm eff}$ is the solid curve and the one with `standard' matrix size
$\Nstd$ is the dotted curve.  The distribution of the lowest eigenvalue is
related to the solution of a non-linear ordinary differential equation of
Painlev\'e~VI (see~\cite{ForrWitte04}). The curves for $N_{\rm eff}$ and $\Nstd$
in figure~\ref{fig:p6_vs_firstzeros} are computed for these non-integral values
of $N$ by establishing a numerical differential equation solver for
Painlev\'e~VI, which is developed in~\cite{DHKMS2}.

We observe from figure~\ref{fig:p6_vs_firstzeros} that the distribution of the
lowest eigenvalue of SO$(2\Neff)$ mimics the distribution of the lowest zero of
rank-zero curves from $\FEp(X)$ much better both in the bulk and in the tail of
the distribution than SO$(2\Nstd)$. However, near the origins we still find a
large discrepancy. We address this discrepancy in the next section.

Substituting $\Neff$ into (\ref{SO2N}) we obtain:
\begin{equation} \label{SO2Na}
\frac{\pi}{\Neff} R_1\left(\frac{\pi y}{\Neff}\right) = \left(1 +
\frac{\sin(2\pi y)}{2\pi y}\right) - r_1 \frac{1 + \cos(2\pi
y)}{L} - \frac{4r_1^2}{6}\frac{\pi y \sin(2\pi y)}{L^2} +
O\left(\frac{1}{L^3}\right).
\end{equation}
With the substitution $N = \Neff$ the scaled one-level density of
SO$(2N)$ now agrees with our conjectural answer (\ref{S1f}) for
that of $\FEp(X)$ in the leading and next-to-leading order term.
The method of Bogomolny \emph{et al.} indicates how to proceed to
find agreement down to the next order term of order $L^{-2}$.
However, with our available data for twists up to $X =
\text{400,000}$ we are not working with sufficiently large $X$ to
proceed further with their method.  For completeness, the details
are worked out in \cite{kn:huynh09}.

\section{Cut-off value for the excised ensemble}\label{sect:disc}

\subsection{Calculating the cut-off value}\label{sect:discvalue}

We now develop an argument to determine the cut-off value $c$
in~\reff{eq:genericdisc}.  In~\cite{CKRS02} and~\cite{CKRS06}, Conrey, Keating,
Rubinstein and Snaith developed a method using random matrix theory to
conjecture the asymptotic order of the number of $L$-functions in a family such
as $\FEp(X)$ that vanish at the central point; alternatively, this is
equivalent, on the Birch and Swinnerton-Dyer conjecture, to the asymptotic order
of the number of curves of rank 2 or higher in the associated family of elliptic
curves. Numerical tests support the prediction for the order, but the associated
proportionality constant could not be determined.  We will use numerical
findings in those papers and the method introduced there to arrive at the
cut-off value $c$.

\kommentar
{
\footnote{See~\cite{CFKRS05} for evidence for the
validity of conjectures like (\ref{MomentsEllipticVsOrthogonal})
for various families of $L$-functions.  Random matrix theory
appears to predict the general structure of these moments, but
misses the arithmetic component.  The contribution from arithmetic
can be incorporated using heuristic methods described in
\cite{CFKRS05}.}
}
We begin by reviewing the method of~\cite{CKRS02}, modifying it to meet
our current purpose. Adjusting slightly the notation of
\cite{CKRS02} to ours, let
\begin{equation}
M_E(X, s) = \frac{1}{X^*} \sum_{{0<d \leq X}\atop{L_E(s,\chi_d)\in
\FEp(X)}} L_E(1/2, \chi_d)^s,
\end{equation}
where $X^* = \#\{0<d \leq X\mid L_E(s,\chi_d)\in \FEp(X)\}$ is the number of
terms in the sum above. Following the philosophy set out in~\cite{KeatSnaith00a}
and~\cite{KeatSnaith00b} we expect that, for large $X$ and $N \sim \log X$,
\begin{equation} \label{MomentsEllipticVsOrthogonal}
M_E(X, s) \sim a_s(E) M_O(N, s),
\end{equation}
where
\begin{equation}
M_O(N, s) = \int_{SO(2N)} \Lambda_A(1,N)^s \, dA
\end{equation}
and $a_s(E)$ is an arithmetical expression depending on the
Dirichlet coefficients $\lambda(p)$ of the curve $E$:
\begin{equation}
  \begin{split}
    a_s(E) &=
    \Bigg[\prod_p\Big(1-\frac{1}{p}\Big)^{s(s-1)/2}\Bigg]\\
    &\qquad\times \Bigg[ \prod_{p\nmid M} \frac{p}{p+1} \Bigg(\frac{1}{p}+ \frac{1}{2}
    \Big[ \mathcal{L}_p \Big( \frac{1}{p^{1/2}}\Big)^s +
    \mathcal{L}_p\Big(\frac{-1}{p^{1/2}}\Big)^s\Big]
    \Bigg) \Bigg]\\
    &\qquad\times \mathcal{L}_M\Big(\frac{\pm\omega(E)}{M^{1/2}}\Big)^s.
  \end{split}
\end{equation}
This holds for prime conductor $M$, where $\omega(E)$ is the sign of the functional
equation of~$L_E(s)$.  For our purposes, the $\pm$ in the last line above must
actually be $+$, corresponding to twists by \emph{positive} fundamental
discriminants $0<d\leq X$ (the $-$ sign corresponds to twists by negative
fundamental discriminants $-X\leq d<0$).  We define
\begin{equation}
{\mathcal L}_p(z)  =  \sum_{n}^\infty \lambda(p^n)z^n=(1 -
 \lambda(p)z + \psi_M(p)z^2)^{-1},
\end{equation}
with $\psi_M(p)$ given at (\ref{psi}).

 In (\ref{MomentsEllipticVsOrthogonal}),
$M_O(N, s)$ denotes the moment generating function of the values
$|\Lambda_A(1,N)|$ as $A$ varies in the random matrix ensemble
SO$(2N)$ (the $s$\textsuperscript{th} moment is the expected value
of $|\Lambda_A(1,N)|^s$). For $\RRe(s)>-1/2$, $M_O(N,s)$ can be
explicitly evaluated~\cite{KeatSnaith00b} as
\begin{equation} \label{MomentsOfSO2N}
M_O(N, s) = \int_{SO(2N)} \Lambda_A(1,N)^s \, dA = 2^{2Ns} \prod_{j =
1}^N \frac{\Gamma(N + j-1) \Gamma(s + j -1/2)}{\Gamma(j -
1/2)\Gamma(s + j + N - 1)}.
\end{equation}
Thus, from~\reff{MomentsOfSO2N}, we have information about the value
distribution of the characteristic polynomials. More precisely, we have, for $c
> 0$, that
\begin{equation} \label{PONx}
P_O(N, x) = \frac{1}{2\pi i x} \int_{c-i\infty}^{c+i\infty}
M_O(N,s)x^{-s}ds;
\end{equation}
here $P_O(N,x)$ denotes the probability density for values of the
characteristic polynomials $\Lambda_A(1,N)$ with $A\in$ SO$(2N)$.

For small $x > 0$, the regime of our interest, the major
contribution in the integral on the right side of~\reff{PONx}
comes from the simple pole at $s = -1/2$ in~\reff{MomentsOfSO2N},
thus
\begin{equation} \label{PONxasymptotic}
P_O(N,x) \sim x^{-1/2} h(N)
\end{equation}
where
\begin{equation}
h(N) = \Res_{s=-1/2} M_O(N,s)= 2^{-N} \Gamma(N)^{-1} \prod_{j=1}^N\frac{\Gamma(N + j - 1)
\Gamma(j)}{\Gamma(j - 1/2) \Gamma(j + N -3/2)}.
\end{equation}

We also make use, for large $N$, of the asymptotic
\begin{equation} \label{hNasymptotics}
h(N) \sim 2^{-7/8} G(1/2) \pi^{-1/4} N^{3/8},
\end{equation}
where $G$ is the Barnes $G$-function~\cite{kn:barnes00}.

\newcommand{\PP}{{\mathcal{P}}}
Notice that $P_O(N, x)\,dx$ is the probability that a characteristic
polynomial $\Lambda_A(1,N)$ of $A \in SO(2N)$ takes a value
between $x$ and $x+dx$. Hence
\begin{equation}
{\rm Prob}(0 \le  \Lambda_A(1,N)\leq \rho) = \int_0^\rho P_O(N,x)
dx.
\end{equation}
 With~\reff{PONxasymptotic} we have for small $x$ that
\begin{equation} \label{PONxevaluated}
{\rm Prob}(0 \le \Lambda_A(1,N) \leq \rho) \sim \int_0^\rho
x^{-1/2}h(N) dx = 2 \rho^{1/2}h(N).
\end{equation}

We now define another probability density $P_E(d,x)$.  Due to
\reff{MomentsEllipticVsOrthogonal}, we expect that this is in some
sense a smooth approximation to the probability density  for
elliptic curve $L$-values from $\FEp(X)$ which have fundamental
discriminants around~$d$. We define, in analogy with~\eqref{PONx},
\begin{equation} \label{ansatzPE}
P_E(d,x) := \frac{1}{2 \pi i x} \int_{c-i\infty}^{c + i
\infty}a_s(E) M_O(\log d, s) x^{-s}ds \sim a_{-1/2}(E) P_O(\log
d,x).
\end{equation}
The final approximation above holds for small $x$; it is obtained by
shifting the line of integration left past the pole at $s = -1/2$
and picking up the respective residue of the integrand.

As described in section~\ref{sect:model}, the formula of Waldspurger
\emph{et al.}~\reff{Waldspurgerformula} implies a discretization of
central $L$-values for elements from $\FEp(X)$ given by
\begin{equation} \label{Waldspurgerdiscretisation} L_E(1/2, \chi_d) = 0
  \quad\text{whenever}\quad L_E(1/2, \chi_d) < \frac{\kappa_E}{\sqrt{d}}.
\end{equation}
By calculating the probability that a random variable, $Y_d$, with probability
density $P_E(d,x)$, takes a value less than $\kappa_E/\sqrt{d}$ and summing over
$d$ up to $X$, the authors in~\cite{CKRS02} and~\cite{CKRS06} predicted the
correct order of magnitude for the number of $L$-functions that vanish at the
central point (see Conjecture 5.1 of~\cite{CKRS06}). However, the constant
factor could not be predicted correctly to agree with the numerical data. Since
their cut-off value $\kappa_E/\sqrt{d}$ did not give the correct number of
vanishing $L$-functions, we will presently work backwards and use the
numerically calculated number of $L$-functions that are zero at the central
point (numerical data from~\cite{CKRS06}) to deduce an ``effective" cut-off
$\tfrac{\delta\cdot\kappa_E}{\sqrt{d}}$ that {\em does} give the correct answer,
for some $\delta>0$. Thus we repeat the calculation of~\cite{CKRS02} but with
the modified cut-off value.

We have
\begin{equation}
  \begin{split}
{\rm Prob}\left(0 \le Y_d \leq \frac{\delta \cdot
    \kappa_E}{\sqrt{d}}\right) &\sim \int_0^{\delta \kappa_E
  d^{-1/2}}a_{-1/2}(E)x^{-1/2}h(\log d) dx\\
&= 2 a_{-1/2}(E) \left(\frac{\delta \cdot \kappa_E}{\sqrt{d}}\right)^{1/2}h(\log
d).
\end{split}
\end{equation}

Now we follow~\cite{CKRS02} and conjecture that
\begin{multline}
  \# \{L_E(s,\chi_d) \in \FEp(X), d {\rm \;prime} : L_E(1/2,\chi_d) = 0 \} =
  \sideset{}{^*}{\sum}_{d\leq X \atop d\ {\rm prime}} {\rm Prob}\left(0 \le Y_d
    \leq \frac{\delta \kappa_E}{\sqrt{d}}\right) \\
  \sim \frac{1}{4\log X} \sum_{n = 1}^{\lfloor X\rfloor} 2a_{-1/2}(E)
  \left(\frac{\delta \kappa_E}{\sqrt{n}}\right)^{1/2} h(\log X),
\end{multline}
where the starred sum means that $d$ is restricted to {\em prime} fundamental
discriminants for which $\chi_d(-M)\omega(E)=+1$ (of which, asymptotically,
there are $X/(4\log X)$ of size at most $X$)%
. Using~\reff{hNasymptotics}, we get
\begin{multline} \label{prelimaryconjecture5_1}
  \# \{L_E(s,\chi_d)\in \FEp(X), d {\rm \;prime} : L_E(1/2,\chi_d) = 0 \} \\
  \sim \frac{1}{4\log X}\cdot 2\; a_{-1/2}(E)
  \sqrt{\kappa_E}\;2^{-7/8}G(1/2)\pi^{-1/4} (\log X)^{3/8} \delta^{1/2}\cdot
  \tfrac{4}{3} X^{3/4}.
\end{multline}
We wish to obtain a numerical value for $\delta$. We would like to thank Michael
Rubinstein for sharing his data from~\cite{CKRS06} on the number of
$L$-functions that vanish at the central point.  For a large number of elliptic
curve families he computes the left side of~\reff{prelimaryconjecture5_1} and
divides it by
\begin{equation}
\frac{1}{4}a_{-1/2}(E)\sqrt{\kappa_E}X^{3/4}(\log X)^{-5/8}.
\end{equation}
The results are plotted (see~\cite{CKRS06}) as a function of $X$; the curves
flatten out and seem to approach a constant value.  One such constant is:
\begin{equation}
0.2834620 {\rm \;\; for\;\;} E=11A_r
\end{equation}
where the nomenclature of the elliptic curve being twisted to
form the family refers to Table 3 of~\cite{CKRS06}.  We will refer to this family as twists of $E_{11}$.

Thus for twists of $E_{11}$ we have
\begin{equation}
\tfrac{8}{3} 2^{-7/8}G(1/2)\pi^{-1/4} \delta^{1/2} \approx 0.2834620. 
\end{equation}
With $G(1/2)$ evaluated as approximately 0.603244, we have
\begin{equation} \label{numericalvaluefordelta}
\delta \approx  0.185116.
\end{equation}
For this family we have
\begin{equation} \label{numericalvalueforkappaE}
\kappa_E = 6.346046521 \quad\mbox{and}\quad a_{-1/2}(E) = 0.732728078.
\end{equation}
Thus, with $\delta$ given in~\reff{numericalvaluefordelta} and
$\kappa_E$ in~\reff{numericalvalueforkappaE}, we take the
`effective' cut-off to be
\begin{equation} \label{delta_in_L_E11}
\delta \cdot \kappa_E = 1.17475.
\end{equation}
That is, when modeling the distribution of $L$-values with
$P_E(d,x)$, as described in~\cite{CKRS02}, integrating up to the value of $\delta
\kappa_E/\sqrt{d}$ gives the correct number of $L$-functions
taking the value zero at the central point.

\subsection{Numerical evidence for the cut-off value: Standard $\Nstd$}\label{sect:discevidenceN0}

The probability densities $P_E$ and $P_O$ are related, for small
values of $x$ (see (\ref{PONxasymptotic}) and (\ref{ansatzPE})),
by
\begin{equation} \label{relation_x_O_vs_x_E}
P_E(d,x) \sim P_O\left(\Nstd, a_{-1/2}^{-2}(E)x\right)
\end{equation}
when $\Nstd\sim\log d$.  Thus a cut-off of $\delta
\kappa_E/\sqrt{d}$ applied to $P_E(d,x)$ scales to
\begin{equation}
c_{{\rm std}}\times \exp(-\Nstd/2): = a_{-1/2}^{-2}(E) \; \delta \; \kappa_E
\times \exp(-\Nstd/2)
\end{equation}
for the distribution of values of characteristic polynomials of
matrices of size $\Nstd$.
Substituting the numerical values we obtain
\begin{equation}\label{eq:cprime}
c_{{\rm std}} \approx 2.188.
\end{equation}

We now present data using the standard matrix size $\Nstd$ to
model the zero statistics. Although  we expect that the excised
orthogonal ensemble models $L$-functions with discriminant around
the value~$X$, for numerical tests we take \emph{all}~$d\leq X$
and set $\Nstd\sim\log X$ so that we have a substantial data set.

\begin{figure}[h]
\includegraphics[scale=0.65]{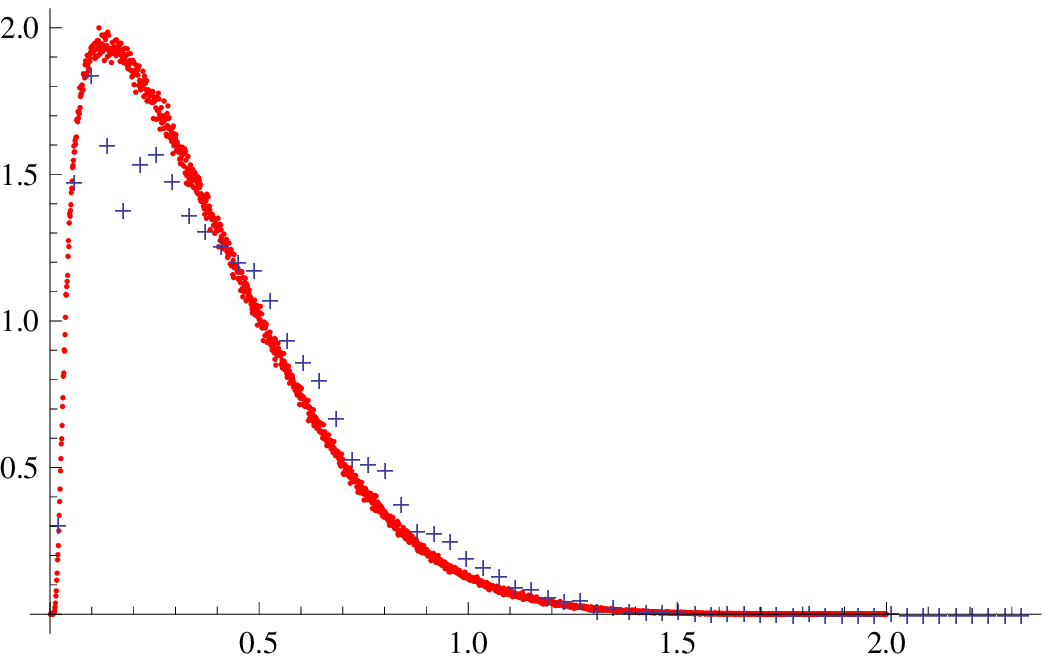}\hspace{0.1
in}
\includegraphics[scale=0.65]{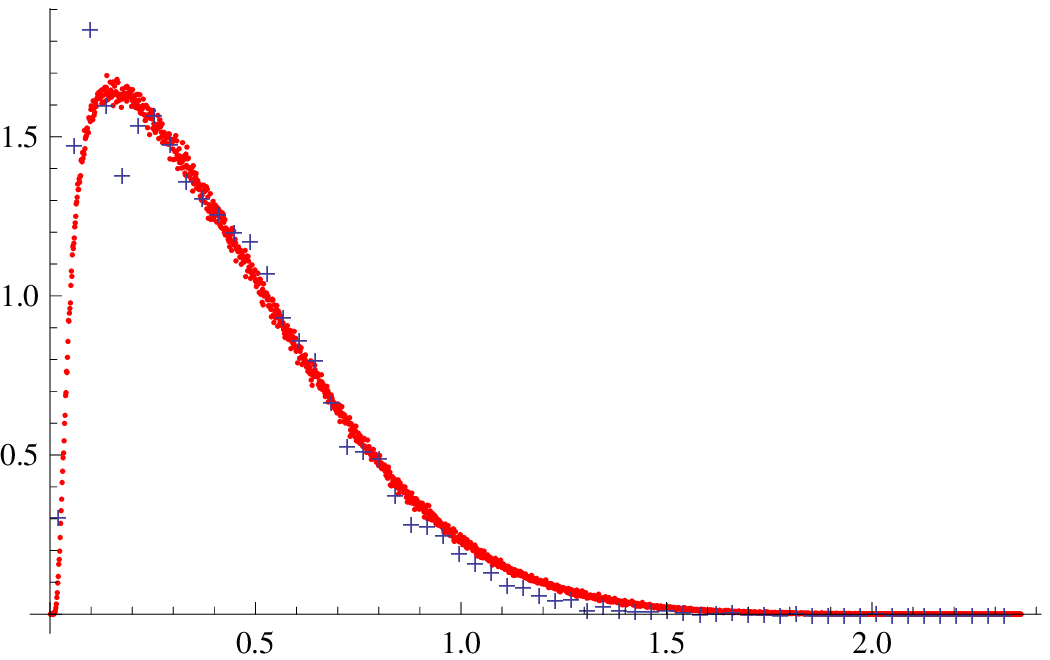}
\caption{Probability density  of the first eigenvalue from $3
\times 10^6$ numerically generated matrices $A \in SO(2\Nstd)$ with
$|\Lambda_A(1,\Nstd)| \geq 2.188\times \exp(- \Nstd/2)$ and $\Nstd=12$ (red
dots) compared with the first zero of even quadratic twists $L_{E_{11}}(s,
\chi_d)$  with  prime fundamental
discriminants $0 < d \leq \text{400,000}$ (blue crosses). In the left
picture the random matrix data is not scaled, in the picture on
the right the mean of the distribution is scaled to match that of
the zero data.} \label{fig:mean}
\end{figure}

\begin{figure}[h]
\includegraphics[scale=1]{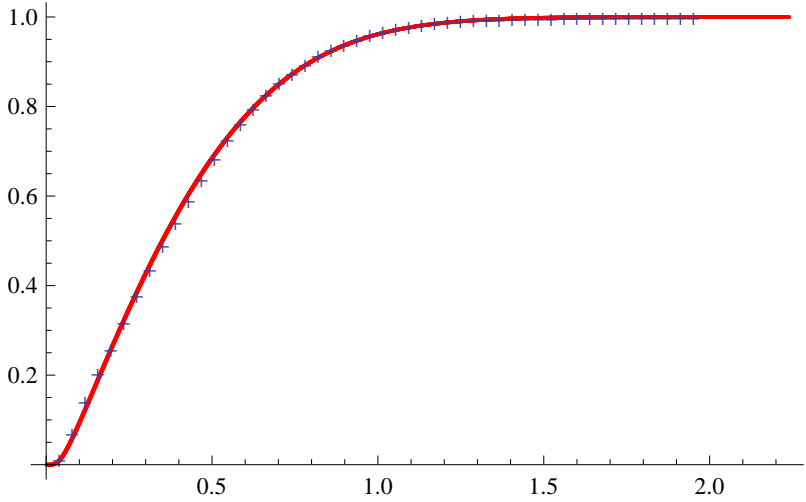}
\caption{Cumulative probability density  of the first eigenvalue
from $3 \times 10^6$ numerically generated matrices $A \in SO(2\Nstd)$
with $|\Lambda_A(1,\Nstd)| \geq 2.188\times \exp(- \Nstd/2)$ and $\Nstd=12$ (red
dots) compared with the first zero of even quadratic twists $L_{E_{11}}(s,
\chi_d)$ with  prime fundamental
discriminants $0 < d \leq \text{400,000}$ (blue crosses). The random
matrix data is scaled so that the means of the two distributions
agree.} \label{fig:cdfmean}
\end{figure}

As mentioned in section~\ref{sect:neff}, and illustrated in
figure~\ref{fig:p6_vs_firstzeros}, eigenvalues of $SO(2N)$
matrices with $N\sim\log X$ do not give particularly good
agreement with the statistics of zeros of $L$-functions from
$\FEp(X)$.  The procedure of excising matrices with small values
of $|\Lambda_A(1,N)|$ from the ensemble does not substantially
change the bulk of the distribution. However, for $E_{11}$, we
illustrate in figure ~\ref{fig:mean} how we can scale the mean of
the distribution of the first eigenvalue to obtain better
agreement. The distribution of the first zero of even quadratic
twists of $L_{E_{11}}(s)$ by prime discriminants between 0 and
400,000 has a mean of 0.4081. The distribution of the first
eigenvalue of $3\times 10^6$ matrices from SO$(2\Nstd)$
conditioned to have $\Lambda_A(1,\Nstd)\geq 2.188 \times
\exp(-\Nstd/2)$, with $\Nstd=12$ ($\approx
\log(\sqrt{11}\;\text{400,000}/2\pi)$) was found to be 0.365;
deviation from the distribution of the first critical zero of the
quadratic twists is quite visible, as illustrated in the left
graphic of figure~\ref{fig:mean}. Rescaling the random-matrix mean
to match that of the $L$-function zeros gives the graphic on the
right side of figure~\ref{fig:mean}, showing much better
agreement.  The cumulative plot of the scaled distribution is
shown in figure~\ref{fig:cdfmean}.

\begin{figure}[h]
\includegraphics[scale=0.9]{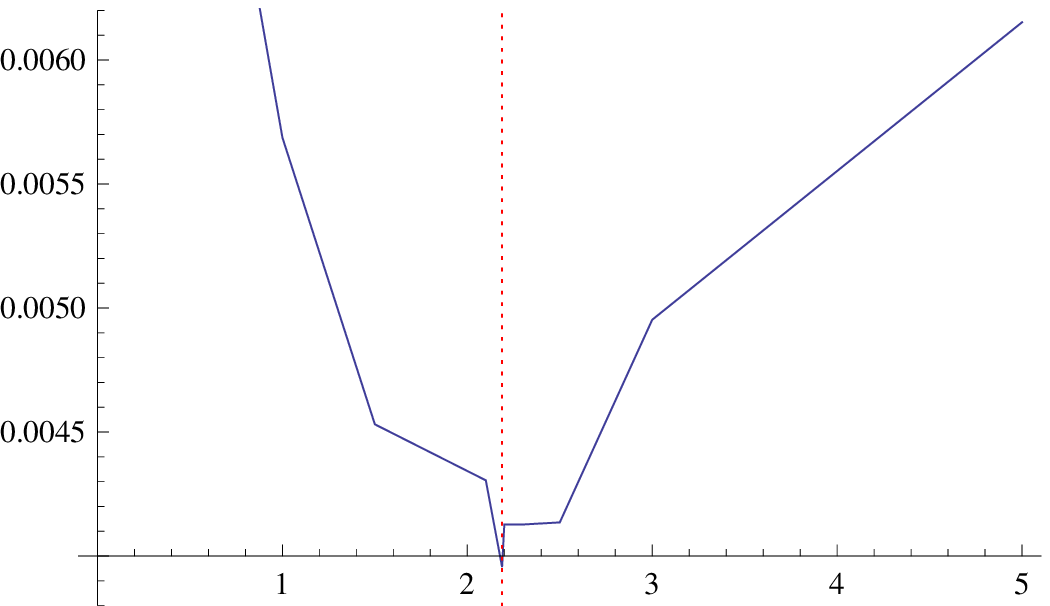}
\caption{A measure of the difference between the mean-scaled
  cumulative distribution of the first eigenvalue of the excised
  random matrix model (\ref{eq:cprime}) with $\Nstd=12$ for various
  values of $c$ and the cumulative distribution of the first zero of
  even quadratic twists $L_{E_{11}}(s, \chi_d)$ with prime fundamental
  discriminants $0 < d \leq \text{400,000}$.  The values $c_{{\rm
      std}}=2.188$ is marked with a vertical line. }
\label{fig:errorE11N}
\end{figure}

Working with the mean-scaled RMT results for $N=\Nstd$, we test various values
of $c$ in (\ref{eq:genericdisc}).
We measure the agreement with the zero distribution by averaging
the absolute value of the difference between the cumulative
distribution of the zeros and the cumulative distribution of the
eigenvalues at a set of evenly spaced points (the blue crosses on
figure~\ref{fig:cdfmean}).  The plot of the difference between the
cumulative distributions versus the cut-off parameter is
shown in figure~\ref{fig:errorE11N}, where $c$ varies along the
horizontal axis. We see a minimum at $c_{{\rm std}}=2.188$; this is the value
predicted in section~\ref{sect:discvalue} and is marked with a dotted
vertical line.

Initially we thought that an
equally plausible way to calculate a cut-off value would be
to equate not the probability densities, as in
(\ref{relation_x_O_vs_x_E}), but rather
the probability of finding an $L$-value of zero with the
probability of discarding a random matrix on the basis of the
condition (\ref{eq:genericdisc}).  Equating probabilities
would lead to a calculation such as
\begin{equation}\label{eq:cprimeprob}
  \begin{split}
    \int_{0}^{\delta \kappa_E/\sqrt{d}} P_E(d,x)dx &\sim \int_0^{\delta
      \kappa_E/\sqrt{d}} a_{-1/2}(E) x^{-1/2} h(\log d)\,dx \\
    &=\int_0^{a_{-1/2}^2(E)\delta \kappa_E /\sqrt{d}} y^{-1/2} h(\Nstd)\,dy \\
    &\sim \int_0^{a_{-1/2}^2(E)\delta \kappa_E /\sqrt{d}} P_O(\Nstd,y)\,dy,
  \end{split}
\end{equation}
implying a cut-off of
\begin{equation}
    a_{-1/2}^2(E)\delta \kappa_E \exp(-\Nstd/2) \approx 0.6307 \exp(-\Nstd/2).
\end{equation}
From figure ~\ref{fig:errorE11N}  we see that $c=0.6307$ certainly
does not minimize the error between the excised random matrix model
and the zero statistics.

It is noteworthy that matching probability densities, rather than
probabilities, appears to be the correct model. Taking the cut-off
value $c_{{\rm std}}$ that best models the zero data (as
calculated at the beginning of this section by matching densities)
means that the proportion of the matrices that are being excluded
by the condition (\ref{eq:genericdisc}) is not the same as the
proportion of $L$-functions from our family of quadratic twists
that are excluded from the zero statistics because they have rank
higher than zero.  The reason for this, and whether the excluded
matrices might potentially model the $L$-functions associated with
the higher-rank curves, are topics for future investigation.


\subsection{Numerical evidence for the cut-off value: $\Neff$}\label{sect:discevidenceNeff}

The excised model can be applied using orthogonal matrices of any
even size~$2N$. We will presently use a method identical to that
in the previous section; however, now we find the cut-off for
matrices of size $N=\Neff$ (which was calculated in
section~\ref{sect:neff}). We recall the shape of $P_O(N,x)$ for
small $x$ and large $N$, equations (\ref{PONxasymptotic}) and
(\ref{hNasymptotics}), and see
\begin{equation}\label{eq:Neffdisc}
  \begin{split}
    P_E(d,x) \sim a_{-1/2} x^{-1/2} h(\log d)&\sim a_{-1/2} x^{-1/2}
    h(2r_1 \Neff)\\
    &\sim a_{-1/2} x^{-1/2} (2r_1)^{3/8} h(\Neff) \\
    & \sim P_O\left( \Neff,a_{-1/2}^{-2}(E)\;(2r_1)^{-3/4}x\right).
  \end{split}
\end{equation}

Thus a cut-off of $\delta \kappa_E/\sqrt{d}$ applied to
$P_E(d,x)$ scales to
\begin{equation}
c_{\rm{eff}}\times \exp(-\Nstd/2): = a_{-1/2}^{-2}(E)(2r_1)^{-3/4} \;
\delta \; \kappa_E \times \exp(-r_1 \Neff)
\end{equation}
for the distribution of values of characteristic polynomials of
matrices of size~$\Neff$.
Plugging in the numerical values we obtain, with $r_1=2.8600$ for
$E_{11}$ (given at (\ref{eq:r1E11}))
\begin{equation}\label{eq:c}
c_{\rm{eff}} \approx 0.5916.
\end{equation}
For this matrix size we use the cut-off value $c=c_{\rm{eff}}$ in
\reff{eq:genericdisc}.

We now compare data for zeros of elliptic curve $L$-functions to
matrices of size~$\Neff$.

We note here that, although the most accurate description of our
expectation is that the excised orthogonal ensemble models
$L$-functions with discriminant around the value $d$ when
$\Neff\sim\log d/(2 r_1)$, for numerical tests we take all $d$ in
an interval such as $0<d\leq X$ and set $\Neff\sim\log X/(2 r_1)$.
This is the most effective way to acquire enough data to have good
resolution in the plots.

\begin{figure}[h]
\includegraphics[scale=0.65]{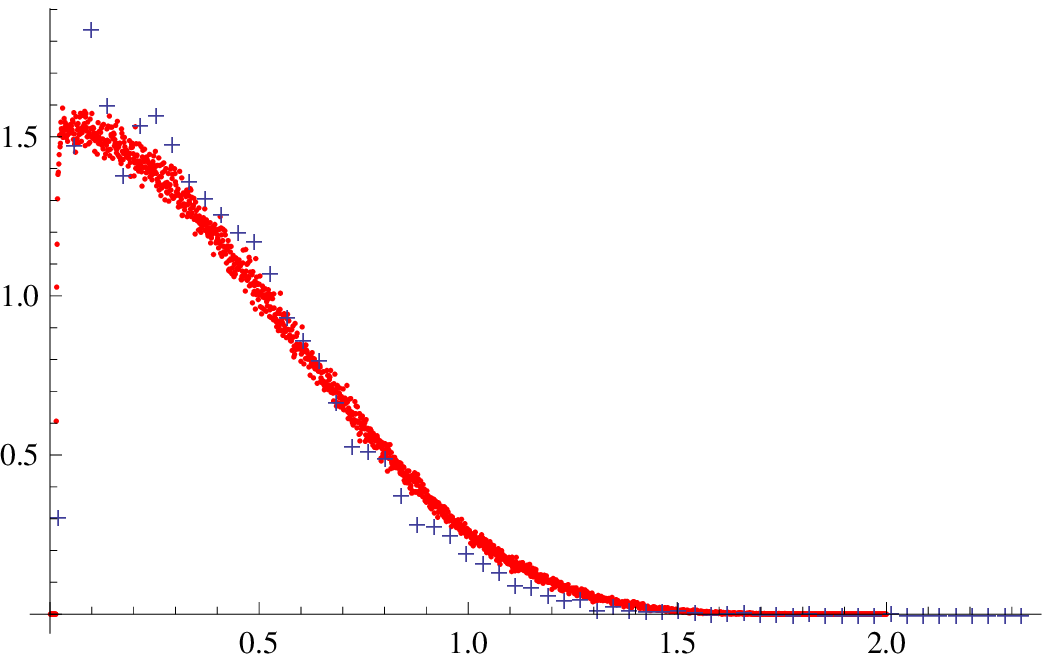}\hspace{0.1
in}
\includegraphics[scale=0.65]{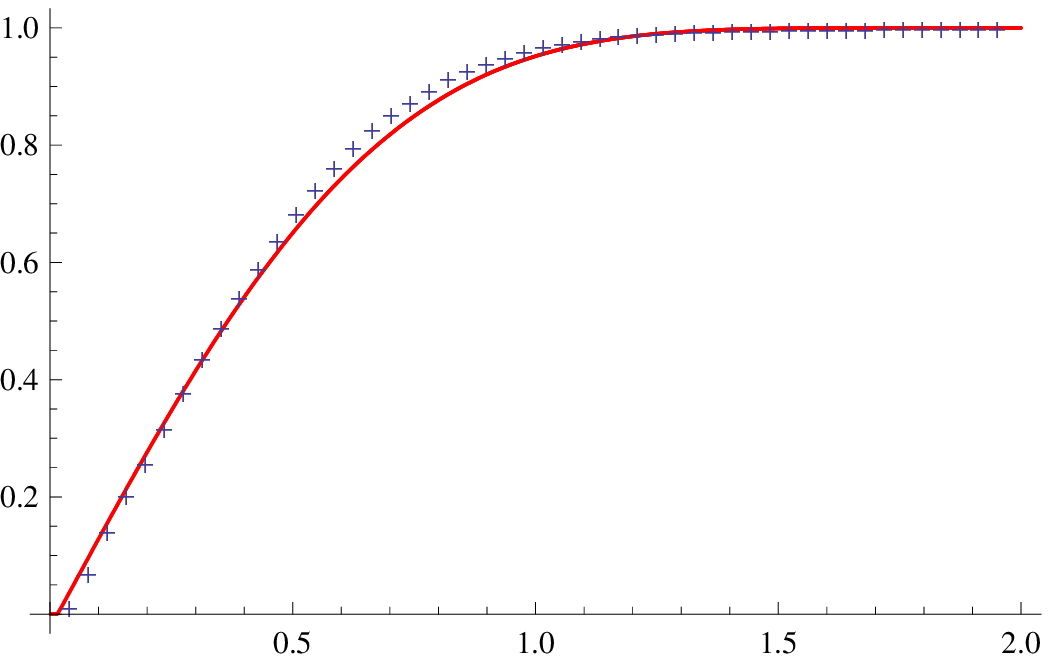}
\caption{Probability density (left) and cumulative probability density
  (right) of the first eigenvalue from $3 \times 10^6$ numerically
  generated matrices $A \in SO(2\Neff)$ with $|\Lambda_A(1,\Neff)|
  \geq 0.5916\times \exp(- r_1\Neff)$ and $\Neff=2$ (red dots)
  compared with the first zero of rank zero even quadratic twists
  $L_{E_{11}}(s,\chi_d)$ with prime fundamental discriminants $0 < d
  \leq \text{400,000}$ (blue crosses). } \label{fig:E11Neff0o5916}
\end{figure}

For $E_{11}$ we have computed zeros of \emph{even} twists
$L_{E_{11}}(s,\chi_d)$ for prime $d$ between 0 and 400,000 which
do \emph{not} vanish at the central point.  We compare the
distribution of the first zero above the critical point of these
$L$-functions with the numerically generated distribution of first
eigenvalues of $3 \times 10^6$ matrices from the excised ensemble
SO(4) (in section~\ref{sect:neff} we computed $\Neff=2.14$) for
various values of the cut-off $c$. Figure~\ref{fig:E11Neff0o5916}
shows the probability density and the cumulative probability
density of first zeros and eigenvalues for $c_{\rm{eff}}=0.5916$.
In contrast to modeling with matrices of size $\Nstd$, here we do
{\em not} scale the mean of the distribution; however, for those
interested, we note that the probability density of the zeros has
mean 0.4081 and the probability density of the eigenvalues has
mean 0.4234.

\begin{figure}[h]
\includegraphics[scale=0.9]{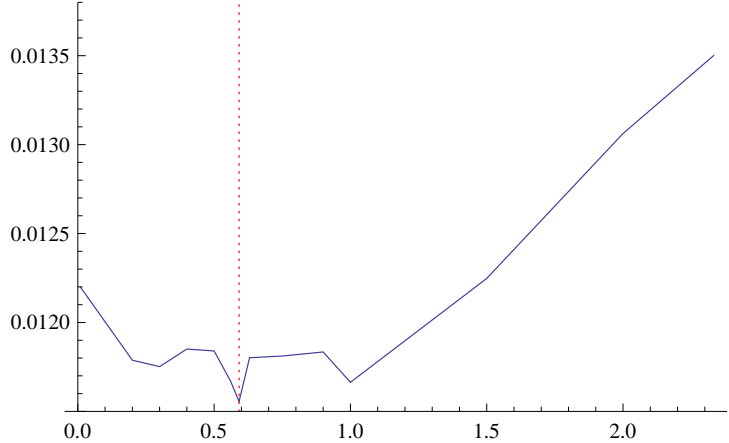}
\caption{A measure of the difference between the cumulative
distribution of the first eigenvalue of the excised random matrix
model (\ref{eq:genericdisc}), as $c$ varies along the horizontal
axis, and the cumulative distribution of the first zero of rank
zero even quadratic twists $L_{E_{11}}(s,\chi_d)$ of with prime
fundamental discriminants $0 < d \leq \text{400,000}$. The value
$c_{\rm{eff}}=0.5916$ is marked with a vertical line. }
\label{fig:errorE11Neff}
\end{figure}

In figure~\ref{fig:errorE11Neff} we plot a measure of the
difference between the first zero distribution and the
distribution of the first eigenvalue for various values of the
cut-off $c$.  Here the error is calculated by
summing the absolute value of the difference between the
cumulative distribution of the zeros and the cumulative
distribution of the eigenvalues at a set of evenly spaced points
(the positions of the blue crosses in figure
\ref{fig:E11Neff0o5916}) and dividing by the number of points.  We
see a minimum at $c_{\rm{eff}}=0.5916$; this is the value predicted in
section~\ref{sect:discvalue} and is marked with a dotted vertical line.

As in the previous section, another possible model is to equate probabilities rather than probability densities. For $\Neff$ this leads to
\begin{equation}
c=(2r_1)^{3/4} a_{-1/2}^{2}(E) \delta \kappa_E\approx 2.3328.
\end{equation}
Again, from figure~\ref{fig:errorE11Neff} we see that this does not come close to
minimizing the error.


\subsection{ Excised model and the one-level density}\label{sect:discdiscussion}
In section~\ref{sect:neff} we provided evidence that the lower-order
terms of the one-level density of $\FEp(X)$ determine an effective
matrix size $\Neff$ such that the distribution of first eigenvalues of
SO(2$\Neff$) models the bulk and tail of the distribution of first
zeros of $\FEp(X)$. In sections
\ref{sect:discvalue}--\ref{sect:discevidenceNeff} we discussed how to
obtain an appropriate cut-off value for the characteristic polynomial
$\Lambda_A$ at 1 with $A \in$ SO$(2N)$ so that the distribution of
first eigenvalues of this subset of SO$(2N)$ models the region at and
near the origin of the distribution of first zeros of $\FEp(X)$. Here
the matrix size is either the effective one, i.e., $N = \Neff$ or the
standard one, i.e., $N = \Nstd$. In the first case no further scaling
is undertaken whereas in the latter one we match mean densities of
SO(2$\Nstd$) to $\FEp(X)$ to achieve qualitative and quantitative
agreement throughout the origin, bulk and tail between the cumulative
distribution of first eigenvalues of SO$(2N)$ and cumulative
distribution of first zeros of $\FEp(X)$. In this section we explore
to what extent the values for the matrix size and the cut-off value on
the characteristic polynomial used for the distribution of first
eigenvalues (to model the distribution of first zeros) can be applied
to the one-level density of eigenvalues (to model the one-level
density of zeros).  We compare zero data of the family $\FEp(X)$ of
quadratic twists $L_{E_{11}}(s, \chi_d)$ with prime fundamental
discriminants $0 < d < X = 400,000$ to eigenvalue statistics of
random matrices from SO(2$N$) for two values of $N$.
\begin{figure}[h]
\begin{center}
\includegraphics[scale=1]{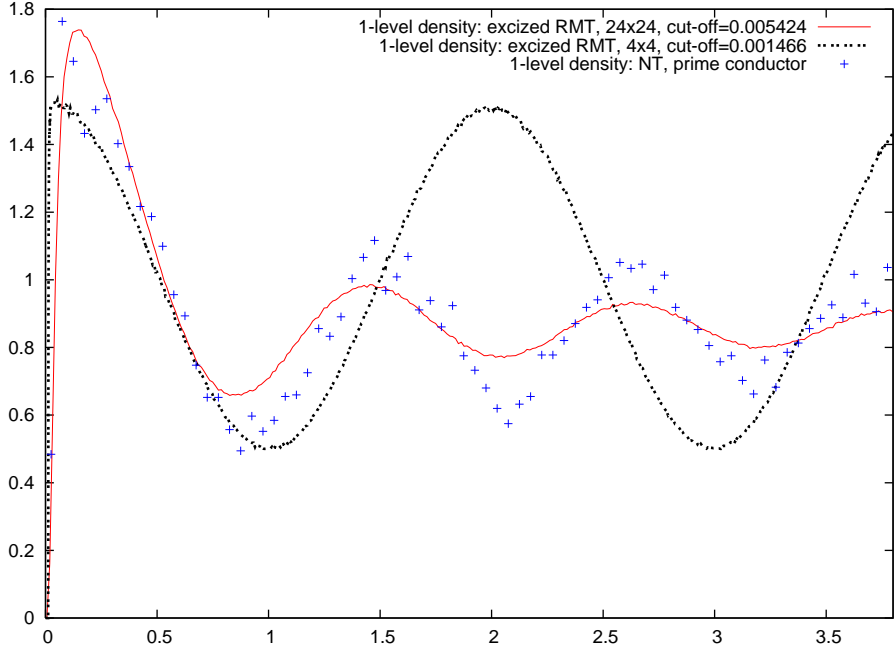}
\caption{\emph{Red solid line:} One-level density of $9.12 \times
  10^6$ numerically generated matrices from SO(2$\Nstd$) and cut-off
  $|\Lambda_A(1,\Nstd)|\geq 2.188 \times \exp(-\Nstd/2) = 0.005424$,
  scaled so that the mean of the first eigenvalue matches that of the
  first zero. \emph{Black dotted line:} One-level density of $9.8
  \times 10^6$ numerically generated matrices from SO(2$\Neff$) and
  cut-off $|\Lambda_A(1,\Neff)|\geq 0.5916 \times \exp(-\Nstd/2) =
  0.001466$. \emph{Blue crosses:} One-level density for zeros of
  even quadratic twists $L_{E_{11}}(s,\chi_d)$ with prime fundamental
  discriminant between 0 and 400,000.}
\label{fig:OLD_scaled_standard_effective_NT_data}
\end{center}
\end{figure}
Figure~\ref{fig:OLD_scaled_standard_effective_NT_data} depicts as
a red solid line the one-level density of matrices $A$ from
SO(2$\Nstd$) with characteristic polynomials constrained to obey
$|\Lambda_A(1,\Nstd)| \geq 2.188 \times \exp(-\Nstd/2) = 0.005424$
with $\Nstd = 12$.  We obtained this by generating $9.12 \times
10^6$ matrices from SO(24).  This one-level density on the random
matrix theory (RMT) side has been scaled, in the horizontal
direction, with the scaling factor $0.4081/0.365 = 1.118$ which we
obtained numerically when comparing with the distribution of the
first zeros in section \ref{sect:discevidenceN0}.  The one-level
density of matrices $A$ from SO(2$\Neff$) with $\Neff=2$ where the
characteristic polynomials are constrained to obey
$|\Lambda_A(1,\Neff)| \geq 0.5916 \times \exp(-\Nstd/2) =
0.001466$ is depicted as a dotted line. We obtained this by
generating $9.8 \times 10^6$ matrices from SO(4). No scaling of
means has been undertaken in this case. The one-level density for
zeros of rank-zero, even quadratic twists of $E_{11}$ with prime
discriminant between 0 and 400 000 is depicted as blue crosses. In
the $\Nstd$ case we find correlation with the zero data over a
wide range whereas in the $\Neff$ case we only find some agreement
up to the first unit mean spacing. Thereafter we see a clear
discrepancy.

\begin{figure}[h]
\begin{center}
\includegraphics[scale=1]{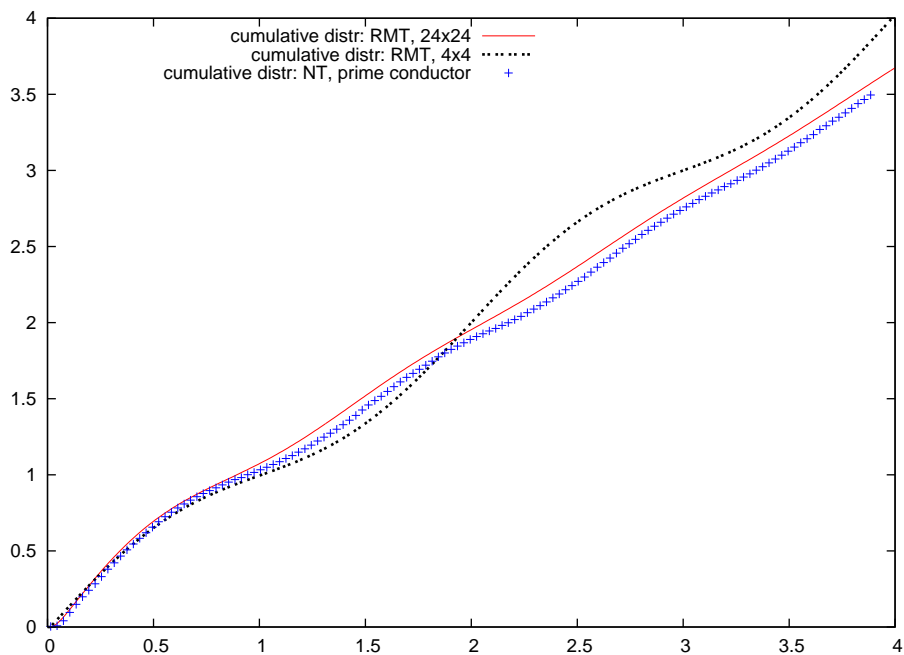}
\caption{Cumulative distributions of one-level densities from figure
 ~\ref{fig:OLD_scaled_standard_effective_NT_data}. \emph{Red solid line:}
  numerically generated matrices from SO(2$\Nstd$), with scaling; \emph{dotted
    line:} numerically generated matrices from SO(2$\Neff$); \emph{blue
    crosses:} zeros of even quadratic twists $L_{E_{11}}(s,\chi_d)$ with prime
  fundamental discriminants between 0 and 400,000.}
\label{fig:OLDstandardtwo}
\end{center}
\end{figure}

\begin{figure}
\begin{center}
\includegraphics[scale=1]{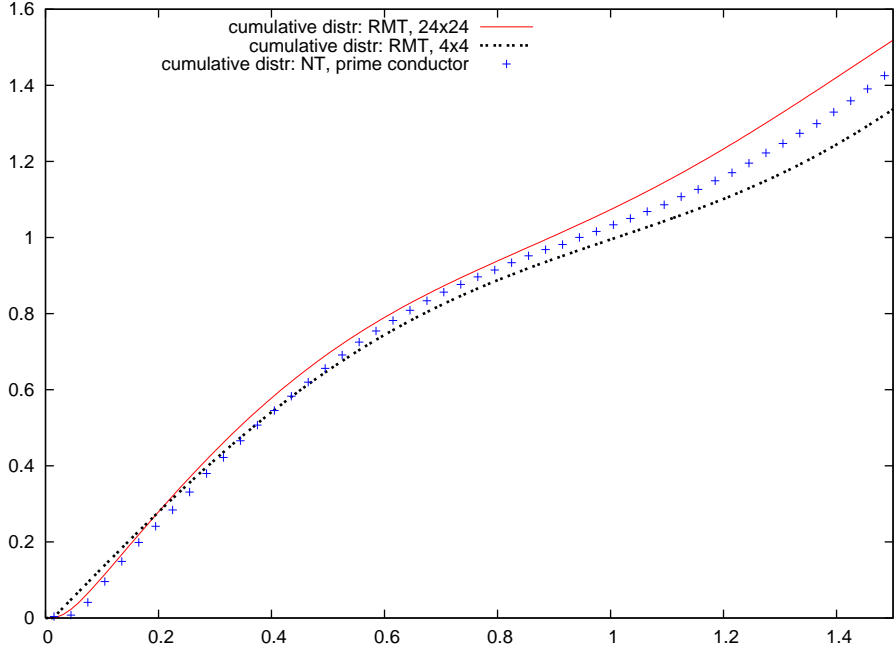}
\caption{Magnification of figure~\ref{fig:OLDstandardtwo}, i.e., cumulative
  distributions of one-level densities from figure
 ~\ref{fig:OLD_scaled_standard_effective_NT_data}. \emph{Red solid line:}
  numerically generated matrices from SO(2$\Nstd$), with scaling; \emph{dotted
    line:} numerically generated matrices from SO(2$\Neff$); \emph{blue
    crosses:} zeros of even quadratic twists $L_{E_{11}}(s,\chi_d)$ with prime
  conductor between 0 and 400,000.}
\label{fig:OLDstandardthree}
\end{center}
\end{figure}

In figure~\ref{fig:OLDstandardtwo} we consider the associated
cumulative distributions from
figure~\ref{fig:OLD_scaled_standard_effective_NT_data}. Here we
observe that the distributions of $\Nstd$ (red line) and the zero
data (blue crosses) track each other over a wide range, whereas
the distributions of $\Neff$ (dotted line) and the zero data (blue
crosses) track each other only up to the first unit spacing.
Thereafter the behaviour of the zeros is not captured.
\newpage Figure \ref{fig:OLDstandardthree} is a magnification of
figure \ref{fig:OLDstandardtwo}. Here we observe that also the
slope of the zero data near the origin is nicely captured by the
excised RMT model with standard matrix size $\Nstd$. This feature
is not so clear in the $\Neff$ case.

In summary, to answer the question to what extent the values for
the matrix size and the cut-off value of the characteristic
polynomial used to model the distribution of first zeros by the
distribution of first eigenvalues  can be applied to the one-level
density, we find that we achieve better agreement when choosing
the excised ensemble with standard matrix size (with the data
scaled so that the mean value of the first eigenvalue matches that
of the first zero) than choosing the excised ensemble with
effective matrix size. In the latter case we only obtain agreement
up to the first unit mean spacing and in the former we get
agreement over a wide range.  It should be noted here that we only
require our model to give useful predictions over a distance of
one mean spacing, because further from the origin we know from
\cite{HKS} that the one level density is strongly dominated by
arithmetic contributions that are accurately modelled by methods
directly incorporating number theoretical information.


\section{One-level density for the excised random matrix
model}\label{sect:onelevel}

In the previous sections we provided evidence that we can use
eigenvalue statistics of random matrices from the excised
orthogonal ensemble to model the zero statistics for the family
$\FEp(X)$ of even quadratic twists of a fixed elliptic curve
$L$-function. Our goal in this section is to compute the one-level
density $R_1^{T_\XX}$ for $T_{\XX}$, the ensemble defined in the
introduction. The main results are contained in the three theorems
set out in the introduction. We prove those theorems in this
section.

\subsection{Proof of Theorem~\ref{lemmaone}} \hspace*{\fill}
\begin{proof}
Consider
\begin{multline} \label{eq:R1TXstart}
R_1^{\TX}(\theta_1) :=\CX \cdot N\int_{0}^\pi \cdots\int_0^\pi
H(\log |\Lambda_A(1,N)|-\XX)\times \\
 \times\prod_{j<k} (\cos \theta_j-\cos \theta_k)^2 d\theta_2\cdots d\theta_N,
\end{multline}
which is the one-level density for the set $\TX$.
Here $H(x)$ denotes the Heaviside function
\begin{equation}
H(x) =
\begin{cases}
1 \mbox{~for~} x > 0\\
0 \mbox{~for~} x < 0,
\end{cases}
\end{equation}
and $\CX$ is a normalization constant (which we discuss later).
Next, we replace in~\reff{eq:R1TXstart} the Heaviside function with its integral representation
\begin{equation}
H(x) = \frac{1}{2 \pi i} \int_{c - i\infty}^{c+i\infty} \frac{\exp(r x)}{r} dr,
\end{equation}
where $c > 0$, and observe that
\begin{align}
  \begin{split}
    \Lambda_A(1,N)
    & = \Lambda_A(\exp(i\theta),N)\bigg|_{\theta = 0} = \prod_{j = 1}^N(1-\exp(i\theta_j)) (1-\exp(-i\theta_j)) \\
    & = 2^N \prod_{j=1}^N (1-\cos\theta_j).
  \end{split}
\label{explicitLambda}
\end{align}
Thus, we have
\begin{equation} \label{Heaviside}
  H(\log|\Lambda_A(1,N)|-\XX) = \frac{1}{2\pi i}
  \int_{c - i\infty}^{c + i\infty} 2^{Nr} \frac{\exp(-r\XX)}{r}
  \prod_{j=1}^N(1 - \cos \theta_j)^r dr.
\end{equation}
Substituting~\reff{Heaviside} into~\reff{eq:R1TXstart} gives
\begin{multline} \label{eq:resultfromlemmaone} R_1^{\TX}(\theta_1) =
  \frac{\CX}{2\pi i} \int_{c-i\infty}^{c+i\infty}2^{Nr}
  \frac{\exp(-r\XX)}{r}\;N \int_0^\pi\cdots\int^\pi_0 \prod_{j=1}^N(1-\cos\theta_j)^r \times \\
  \times \prod_{j<k}(\cos \theta_j-\cos \theta_k)^2d\theta_2 \cdots d\theta_N\,dr.
\end{multline}
Now, observe that in~\reff{eq:resultfromlemmaone} we have
\begin{multline}
  N \int_0^\pi\cdots\int^\pi_0 \prod_{j=1}^N(1-\cos \theta_j)^r \prod_{j<k}(\cos \theta_j-\cos \theta_k)^2d\theta_2 \cdots d\theta_N\\
  = N \int_0^\pi\cdots\int^\pi_0 \prod_{j=1}^N w^{(r-1/2,-1/2)}(\cos \theta_j)\prod_{j<k}(\cos \theta_j-\cos \theta_k)^2d\theta_2 \cdots d\theta_N
\end{multline}
where $w^{(\alpha, \beta)}(\cos \theta) = (1-\cos\theta)^{\alpha+1/2}(1+\cos\theta)^{\beta+1/2}$ is the weight function for the Jacobi ensemble of random matrices~\cite{LGRM}.
We now observe that
\begin{equation}
R_1^{J_N}(\theta_1;\alpha,\beta) = \int_0^\pi\cdots\int^\pi_0 \prod_{j=1}^N w^{(r-1/2,-1/2)}(\cos \theta_j)\prod_{j<k}(\cos \theta_j-\cos \theta_k)^2d\theta_2 \cdots d\theta_N
\end{equation}
is the one-level density for the Jacobi ensemble $J_N$. This completes the proof.
\end{proof}

\subsection{Computation of the normalization constant $\CX$}
Recall that when integrating over SO$(2N)$ with respect to the normalized Haar measure
using Weyl's integration formula~\cite{Weyl} the normalization constant $\CSO$ is determined by
\begin{equation}
1 = \int_{{\rm SO}(2N)} dA
= \CSO \int_{[0,\pi]^N} \prod_{j < k} (\cos\theta_j - \cos\theta_k)^2 d \theta_1 \cdots d\theta_N.
\end{equation}
Selberg's integral formula states that for integral $N$ and complex $r,s$ with $\RRe(r)$, $\RRe(s) > -1/2$
\begin{multline} \label{SelbergChapter3}
 \int_0^\pi d\phi_1 \cdots \int_0^\pi d \phi_N \prod_{l = 1}^N (1 -\cos\phi_l)^r (1 + \cos\phi_l)^s
\prod_{1 \leq j < k \leq N} (\cos\phi_j - \cos\phi_k)^2\\
 = 2^{N(N + r + s - 1)} \times \prod_{j=0}^{N-1}\frac{\Gamma(2 +j)\Gamma(s + 1/2+j)
\Gamma(r + 1/2 + j)}{ \Gamma(s + r + N + j)}.
\end{multline}
Using~\reff{SelbergChapter3} the normalization constant has the explicit form
\begin{equation} \label{explicitCSO}
\CSO = 2^{-N(N-1)} \prod_{j = 0}^{N - 1} \frac{\Gamma(N + j)}{\Gamma(2+j) \Gamma(1/2+j)^2}.
\end{equation}

Likewise, the normalization constant $\CX$ is determined by
\begin{equation} \label{condition_on_CX}
1 = \CX \int_{[0, \pi]^N} H(\log\Lambda_A(1,N) - \XX) \prod_{j < k} (\cos \theta_j - \cos\theta_k)^2 d\theta_1 \cdots d\theta_N.
\end{equation}
Using~\reff{Heaviside}
yields
\begin{multline} \label{CSOCXfraction}
   \frac{\CSO}{\CX}
    = \frac{\CSO}{2 \pi i} \int_{c - i\infty}^{c+i\infty} \int_{[0,\pi]^N} d\theta_1 \cdots d\theta_N d\alpha\times \\
    \times 2^{N\alpha}\frac{\exp(-\alpha \XX)}{\alpha}
    \prod_{j=1}^N(1-\cos\theta_j)^\alpha \prod_{j < k}(\cos\theta_j -
    \cos\theta_k)^2.
\end{multline}
Substituting~\reff{explicitCSO} into~\reff{CSOCXfraction} and applying Selberg's integral formula~\reff{SelbergChapter3}
we obtain
\begin{equation} \label{CSOCX}
\frac{\CSO}{\CX} = \frac{1}{2\pi i} \int_{c-i\infty}^{c+i\infty} \frac{\exp(-\alpha\XX)}{\alpha}2^{2N\alpha}
\prod_{j = 0}^{N-1} \frac{\Gamma(N+j)\Gamma(\alpha + 1/2 +j)}{\Gamma(\alpha + N +j) \Gamma(1/2+j)}d\alpha.
\end{equation}
The evaluation of~\reff{CSOCX} boils down to closing the contour to the left and to computing the
residues associated to the poles at $\alpha=0$ and the negative half integers coming from the $\Gamma(\alpha + 1/2 + j)$-term. The residue at the simple pole $\alpha=0$ is 1 and the residue at the simple pole at $\alpha=-1/2$ is
\begin{equation}
-2\exp(\XX/2) 2^{-N} \frac{\Gamma(N)}{\Gamma(N-1/2)\Gamma(1/2)}
\prod_{j=1}^{N-1} \frac{\Gamma(N+j)\Gamma(j)}{\Gamma(N+j-1/2)\Gamma(1/2+j)}.
\end{equation}
We denote the contribution to~\reff{CSOCX} from residues of higher order poles at $\alpha_k = \frac{-1-2k}{2}$ with $k = 1,2,3,\ldots$ by
\begin{equation} \label{eq:contributionhigherorderpoles}
\sum_{k\geq 1} a_k \exp((k+1/2)\chi).
\end{equation}
Thus, we obtain
\begin{equation} \label{CSOasresidues}
  \begin{split}
    \frac{\CSO}{\CX} & = 1-\frac{\exp(\XX/2)}{2^{N-1}}
    \frac{\Gamma(N)}{\Gamma(N-1/2)\Gamma(1/2)} \prod_{j=1}^{N-1}
    \frac{\Gamma(N+j)\Gamma(j)}{\Gamma(N+j-1/2)\Gamma(1/2+j)}
    \\
    & \ \ \ ~~~+\ \sum_{k\geq 1} a_k \exp((k+1/2)\chi).
  \end{split}
\end{equation}
Notice that when $\XX \rightarrow -\infty$ we have $\CX
\rightarrow \CSO$, as expected.  The regime of interest for us is
when $\XX < 0$. From~\reff{eq:contributionhigherorderpoles} we see
that the contributions of residues of higher poles decrease
exponentially for $\XX < 0$.  The computation of residues of
higher order poles is easily done with a computer algebra system
(we used Mathematica). We find that for reasonable $\XX$ we have
good convergence using about 10 poles. When using only the $K$
rightmost poles the error term is $O(\exp(-c_K \XX))$ with
$\alpha_{K+1}< c_K < \alpha_{K}$. We refer to
section~\ref{subsection:proofofcorollaryone} for the relevant
details. There we will also discuss the convergence of a series
similar to~\reff{eq:contributionhigherorderpoles}.
We are then left to show that the contribution from closing the contour
is~0. This is indeed the case for $\exp(\XX) < 2^{2N} = \max[{\Lambda_A(1,N)}]$
but we skip the calculation here and instead also refer to
section~\ref{subsection:proofofcorollaryone} as we do a similar computation for
$\RTX$ there.

\subsection{Proof of Theorem~\ref{theoremone}}
\begin{proof}
  We rewrite~\reff{eq:resultfromlemmaone} so we can apply standard methods in
  the theory of orthogonal polynomials. First we consider the general Jacobi
  ensemble with weight function $(1-x)^\alpha(1+x)^\beta$ and then specialize to
  our setting with $\alpha=r-1/2$ and $\beta=-1/2$. Notice that the weight
  function here differs slightly from the one in section~\ref{introduction}.

Following Szeg\H{o}~\cite{Szego39} we define the Jacobi polynomials
$P_n^{(\alpha,\beta)}(x)$, $\alpha, \beta>-1$, which satisfy
\begin{equation}
\int_{-1}^1 \label{JacobiPolynomial}
P_n^{(\alpha,\beta)}(x)P_m^{(\alpha,\beta)}(x)(1-x)^\alpha
(1+x)^\beta dx = \delta_{nm} h_n^{(\alpha,\beta)},
\end{equation}
where
\begin{equation}
h_n^{(\alpha,\beta)} = \frac{2^{\alpha+\beta+1}} {2n+\alpha +\beta
+1} \frac{\Gamma(n+\alpha+1)\Gamma(n+\beta+1)}
{\Gamma(n+1)\Gamma(n+\alpha+\beta+1)}.
\end{equation}
Changing to angular variables,~\reff{JacobiPolynomial} takes the equivalent form
\begin{equation}
\int_0^\pi P_n^{(\alpha,\beta)}(\cos
\theta)P_m^{(\alpha,\beta)}(\cos \theta)(1-\cos
\theta)^{\alpha+1/2} (1+\cos \theta)^{\beta +1/2} d\theta =
\delta_{nm} h_n^{(\alpha,\beta)}.
\end{equation}

Using the Vandermonde determinant and further matrix row and column operations we have
\begin{multline}
  \prod_{j<k} (\cos \theta_k-\cos \theta_j) = \left|
    \begin{array}{ccccc} 1&\cos \theta_1& \cos^2\theta_1&\ldots &
      \cos^{N-1}\theta_1\\ 1& \cos \theta_2& \cos^2\theta_2&\ldots
      &\cos^{N-1}\theta_2\\\vdots& \vdots&\vdots&\ddots&\vdots\\
      1&\cos\theta_N
      &\cos^2\theta_N&\ldots&\cos^{N-1}\theta_N\end{array}\right|  \\
   \shoveleft{=\prod_{j=0}^{N-1} \frac{(h_j^{(\alpha,\beta)})^{1/2}}
  {\ell_j^{(\alpha,\beta)}}  \times} \\
   \times \left|\begin{array}{cccc} \frac{1}{\sqrt{h_0^{(\alpha,\beta)}}}
      P_0^{(\alpha,\beta)}(\cos\theta_1)&
      \frac{1}{\sqrt{h_1^{(\alpha,\beta)}}}
      P_1^{(\alpha,\beta)}(\cos\theta_1)&\ldots&\frac{1}{\sqrt{h_{N-1}^{(\alpha,\beta)}}}
      P_{N-1}^{(\alpha,\beta)}(\cos\theta_1)\\\frac{1}{\sqrt{h_0^{(\alpha,\beta)}}}
      P_0^{(\alpha,\beta)}(\cos\theta_2)&\frac{1}{\sqrt{h_1^{(\alpha,\beta)}}}
      P_1^{(\alpha,\beta)}(\cos\theta_2)&\ldots&\frac{1}{\sqrt{h_{N-1}^{(\alpha,\beta)}}}
      P_{N-1}^{(\alpha,\beta)}(\cos\theta_2)\\\vdots&\vdots&\ddots&\vdots\\
      \frac{1}{\sqrt{h_0^{(\alpha,\beta)}}}
      P_0^{(\alpha,\beta)}(\cos\theta_{N})&\frac{1}{\sqrt{h_1^{(\alpha,\beta)}}}
      P_1^{(\alpha,\beta)}(\cos\theta_{N})&\ldots&\frac{1}{\sqrt{h_{N-1}^{(\alpha,\beta)}}}
      P_{N-1}^{(\alpha,\beta)}(\cos\theta_N)\end{array} \right|,
\end{multline}
where $\ell_j^{(\alpha,\beta)}$  is the leading coefficient of the
polynomial $P_j^{(\alpha,\beta)}$:
\begin{equation}
\ell_j^{(\alpha,\beta)} =2^{-j}\binom{2j+\alpha+\beta}{j}.
\end{equation}
Thus we have the following determinantal expression
\begin{multline}
  \prod_{j<k} (\cos \theta_k-\cos \theta_j)^2\\
  = \prod_{j=0}^{N-1} \frac{h_j^{(\alpha,\beta)}}
  {(\ell_j^{(\alpha,\beta)})^2}\det\left(\sum_{n=1}^N
    (h_{n-1}^{(\alpha,\beta)})^{-1} P_{n-1}^{(\alpha,\beta)}(\cos
    \theta_j)P_{n-1}^{(\alpha,\beta)}(\cos \theta_k)\right)_{j,k=1,\ldots,N}.
\end{multline}

Setting
\begin{equation} \label{Csr}
C_{s,r}:=2^{-N(N+r+s-1)}\prod_{j=0}^{N-1} \frac{\Gamma(s+r+N+j)}
{\Gamma(2+j)\Gamma(s+1/2+j)\Gamma(r+1/2+j)},
\end{equation}
the normalized measure is
\begin{equation} \label{normalizedmeasure}
  \begin{split}
    &C_{s,r}\prod_{j=1}^N(1-\cos\theta_j)^r(1+\cos\theta_j)^s \prod_{j<k} (\cos
    \theta_k-\cos \theta_j)^2\\
    & = \frac{1}{N!} \det \Big( \big(\sum_{n=1}^N(h_{n-1}^{(r-1/2,s-1/2)})^{-1}
    P_{n-1}^{(r-1/2,s-1/2)}(\cos
    \theta_j)P_{n-1}^{(r-1/2,s-1/2)}(\cos \theta_k)\big)\times \\
    &\qquad\qquad\qquad
    \times(1+\cos\theta_j)^{s/2}(1-\cos\theta_j)^{r/2}(1+\cos\theta_k)^{s/2}
    (1-\cos \theta_k)^{r/2} \Big)_{j,k=1,\ldots,N} \\
    &=\frac{1}{N!}\det\Big(f_N^{(r-1/2,s-1/2)}(\theta_j,\theta_k)\Big)_{j,k=1,\ldots,N},
  \end{split}
\end{equation}
where $f_N^{(r-1/2,s-1/2)}(\theta_j,\theta_k)$ is implicitly defined as the
expression inside `det$(\cdot)$' in the middle term above.

Observe that this determinantal kernel $f_N$ satisfies the hypotheses of
Gaudin's Lemma (see e.g. Theorem 5.2.1 in~\cite{Mehta91}), namely
\begin{equation}
\int_0^\pi f_N^{(r-1/2,s-1/2)}(\theta,\theta)d\theta =N
\end{equation}
and
\begin{equation}
\int_0^\pi f_N(x,\theta)f_N(\theta,y)d\theta = f_N(x,y).
\end{equation}
By Gaudin's Lemma we then have
\begin{equation} \label{eq:Gaudin}
\int_0^\pi \det(f(\theta_j, \theta_k))_{j,k=1,\ldots,N}d\theta_N = (N - (N-1)) \det(f(\theta_j, \theta_k))_{j,k=1,\ldots,N-1}.
\end{equation}
Applying~\reff{eq:Gaudin} $N - n$ times, together
with~\reff{normalizedmeasure}, gives the following formula for the
$n$-level density:
\begin{multline} \label{GaudinAppliedToRn}
  R_n(\theta_1,\cdots,\theta_n) = \frac{N!}{(N-n)!} \int_0^\pi \cdots
  \int_0^\pi \frac{1}{N!} \det (f_N(\theta_j,\theta_k))_{N\times N}
  d\theta_{n+1}\cdots d\theta_N \\
  = \frac{1}{(N-n)!}(N-(N-1))(N-(N-2))\cdots(N-(n+1-1))
  \det(f_N(\theta_j, \theta_k))_{n\times n} \\
  \shoveleft{= \det(f_N(\theta_j,\theta_k))_{n\times n}.}\\
\end{multline}

So, using~\reff{GaudinAppliedToRn} in~\reff{eq:resultfromlemmaone}, we arrive at
\begin{align} \label{MyRTX}
\RTX(\theta)=\frac{\CX}{2\pi i} \int_{c-i\infty}^{c+i\infty}
\frac{2^{Nr}}{C_{0, r}} \frac{\exp(-r\XX)}{r} f_N^{(r-1/2,-1/2)}(\theta,\theta)\,dr
\end{align}
with $C_{0, r}$ given in~\reff{Csr}.

With the Christoffel-Darboux formula (see equation (4.5.2) in~\cite{Szego39}) we have
\begin{multline} \label{eq:cd}
  f_N^{(r-1/2,s-1/2)}(\theta_j,\theta_k)=\frac{2^{-r-s+1}} {2(N-1) +r+s+1}
  \frac{\Gamma(N+1)\Gamma(N+r+s)} {\Gamma(N+r-1/2)\Gamma(N+s-1/2)} \times \\
  \shoveright{\times
    \frac{(1+\cos\theta_j)^{s/2}(1-\cos\theta_j)^{r/2}(1+\cos\theta_k)^{s/2}
      (1-\cos \theta_k)^{r/2}}{\cos \theta_j-\cos \theta_k}\times} \\
  \times \bigg[
  P_N^{(r-1/2,s-1/2)}(\cos\theta_j)P_{N-1}^{(r-1/2,s-1/2)} (\cos \theta_k) \\
  - P_{N-1}^{(r-1/2,s-1/2)}(\cos\theta_j)P_N^{(r-1/2,s-1/2)}(\cos \theta_k)
  \bigg].
\end{multline}

Thus $f(\theta,\theta)$ in~\reff{MyRTX} reduces, with the Christoffel-Darboux formula (\ref{eq:cd}), to
\begin{multline}
   f_N^{(r-1/2,-1/2)} (\theta,\theta)  = (1-\cos \theta)^r
    \frac{2^{1-r}}{2N+r-1} \frac{\Gamma(N+1)\Gamma(N+r)} {\Gamma(N+r-1/2)\Gamma(N-1/2)}\times \\
    \times\Bigg[ \Big[\frac{d}{d\cos \theta} P_N^{(r-1/2,-1/2)}(\cos
    \theta)\Big] P_{N-1}^{(r-1/2,-1/2)}(\cos \theta) \\
    \shoveright{ {}-P_N^{(r-1/2,-1/2)}(\cos \theta) \frac{d}{d\cos\theta}
    P_{N-1}^{(r-1/2,-1/2)} (\cos \theta)\Bigg]}\\
    =(1-\cos \theta)^r \frac{2^{1-r}}{2N+r-1} \frac{\Gamma(N+1)\Gamma(N+r)}
    {\Gamma(N+r-1/2)\Gamma(N-1/2)} P(N,r,\theta),
\label{fNasPolynomial}
\end{multline}
where we define
\begin{multline}\label{PNrtheta}
  P(N,r,\theta) := \Big[\frac{d}{d\cos \theta} P_N^{(r-1/2,-1/2)}(\cos
  \theta)\Big] P_{N-1}^{(r-1/2,-1/2)}(\cos \theta) \\
  {}- P_N^{(r-1/2,-1/2)}(\cos \theta) \frac{d}{d\cos\theta}
  P_{N-1}^{(r-1/2,-1/2)} (\cos \theta).
\end{multline}

With~\reff{fNasPolynomial} in~\reff{MyRTX} we arrive at
\begin{multline}
   R_1^{\TX}(\theta) = \frac{\CX}{2\pi i} \int_{c-i\infty}^{c+i\infty}
    \frac{\exp(-r\XX)}{r} 2^{N^2+2Nr-N} \times \\
    \times \prod_{j=0}^{N-1} \frac{\Gamma(2+j)\Gamma(1/2+j)\Gamma(r+1/2+j)}
    {\Gamma(r+N+j)} \times \\
    \times (1-\cos \theta)^r \frac{2^{1-r}}{2N+r-1}
    \frac{\Gamma(N+1)\Gamma(N+r)} {\Gamma(N+r-1/2)\Gamma(N-1/2)} P(N,r,\theta)\,dr.
\label{eq:onelevelintegral}
\end{multline}
\end{proof}

\subsection{Proof of Theorem~\ref{corollaryone}}\label{subsection:proofofcorollaryone}
\begin{proof}
We consider~\reff{eq:onelevelintegral} (that is the form of $\RTX$ given in Theorem~\ref{theoremone}).
First, we observe that by~\cite{Szego39}, pages 63--64, the Jacobi polynomial $P_N^{(\alpha,\beta)}(x)$ is a polynomial in
$\alpha, \beta$ and $x$ for arbitrary complex values of $\alpha$ and $\beta$.
Hence $P(N, r, \theta)$ is a polynomial in $r$.

The poles arising from $\Gamma(N+r)$ at negative integers $\leq -N$ are
cancelled by the zeros of the term $1/\Gamma(r+N+j)$ for $j=0$.  The pole at
$r=1-2N$ of $1/(2N+r-1)$ is cancelled by the zero of $1/\Gamma(N+r+j)$ when
$j=N-1$. We now discuss under which conditions we can close the contour to the
left. For this, it is helpful to analyze the integrand of
(\ref{eq:onelevelintegral}) when $|r|$ is large. Our main reference for various
identities and formul\ae\ in the following is~\cite{AbramowitzStegun65}.

Using equation (22.5.42) in~\cite{AbramowitzStegun65} we write the Jacobi polynomials
as
\begin{equation}\label{eq:hypergeom}
P_N^{(r-1/2,-1/2)}(\cos \theta)=\binom{N+r-1/2}{N}
F(-N,N+r;r+1/2;\tfrac{1-\cos \theta}{2}),
\end{equation}
where $F$ is a hypergeometric function.  This representation
allows us to calculate the derivative of the Jacobi polynomial
that appears in $P(N,r,\theta)$ using (15.2.1 in~\cite{AbramowitzStegun65})
\begin{equation}
\frac{d}{dz} F(a,b;c;z) = \frac{ab}{c}F(a+1,b+1;c+1;z).
\end{equation}
Combining the last identity with (\ref{eq:hypergeom}) gives
\begin{multline}
  \frac{d}{d\cos \theta} P_N^{(r-1/2,-1/2)}(\cos \theta)\\
  = -\frac{1}{2}
  \binom{N+r-1/2}{N}\frac{d}{dz}F(-N,N+r;r+1/2;z)\bigg|_{z=\tfrac{1-\cos\theta}{2}} \\
  \shoveleft{ =-\frac{1}{2} \binom{N+r-1/2}{N} \frac{(-N)(N+r)}{r+1/2}
    F(-N+1,N+r+1;r+3/2;\tfrac{1-\cos\theta}{2}).}\\
\label{derivativeofPnrtheta}
\end{multline}

Substituting~\reff{eq:hypergeom} and~\reff{derivativeofPnrtheta} into~\reff{PNrtheta} yields
\begin{equation}
  \begin{split}
    P(N,r,\theta) &=
    -\frac{1}{2(r+1/2)}\binom{N+r-1/2}{N}\binom{N+r-3/2}{N}\times\\
    &\qquad \times\Big[
    (-N)(N+r)F(-N+1,N-1+r;r+1/2;\tfrac{1-\cos\theta}{2}) \times \\
    &\quad\qquad\qquad\qquad\qquad \times
    F(-N+1,N+r+1;r+3/2;\tfrac{1-\cos\theta}{2}) \\
    &\qquad\qquad {}-(-N+1)(N-1+r)
    F(-N,N+r;r+1/2;\tfrac{1-\cos\theta}{2})\times \\
    &\quad\qquad\qquad\qquad\qquad\qquad\qquad\qquad \times
    F(-N+2,N+r;r+3/2;\tfrac{1-\cos\theta}{2})\Big].
  \end{split}
\end{equation}

We want the asymptotics for large $|r|$, and this is easier when
$r$ appears only in one argument of the hypergeometric function.
For this, we use identity 15.3.4 from~\cite{AbramowitzStegun65}:
\begin{equation}
F(a,b;c;z)=(1-z)^{-a} F(a,c-b;c;\tfrac{z}{z-1}),
\end{equation}
and so we obtain
\begin{equation}
  \begin{split}
    P(N,r,\theta)
    &= -\frac{1}{2(r+1/2)}\binom{N+r-1/2}{N}\binom{N+r-3/2}{N} \Big(\frac{1+\cos \theta}{2}\Big)^{2N-2}\times \\
    &\qquad \times\Big[
    (-N)(N+r)F(-N+1,-N+3/2;r+1/2;\tfrac{\cos\theta-1}{\cos\theta +1})\times \\
    &\quad\qquad\qquad\qquad\qquad \times
    F(-N+1,1/2-N;r+3/2;\tfrac{\cos\theta-1}{\cos\theta+1}) \\
    &\qquad\qquad {}-(-N+1)(N-1+r)
    F(-N,1/2-N;r+1/2;\tfrac{\cos\theta-1}{\cos\theta +1})\times \\
    &\quad\qquad\qquad\qquad\qquad\qquad\qquad\qquad\times
    F(-N+2,3/2-N;r+3/2;\tfrac{\cos\theta-1}{\cos\theta+1})\Big].
  \end{split}
\end{equation}

Recalling equation 15.7.1 from~\cite{AbramowitzStegun65},
\begin{eqnarray}
&&F(a,b;c;z)=\sum_{n=0}^m \frac{\Gamma(a+n)\Gamma(b+n)\Gamma(c)}
{\Gamma(a)\Gamma(b)\Gamma(c+n)} \frac{z^n}{n!}+O(|c|^{-m-1}),
\end{eqnarray}
which holds for fixed $a$,$b$ and $z$ and large $|c|$ we
get
\begin{equation} \label{Fasymptotics}
F\left(-N+1,-N+3/2;r+1/2;\tfrac{\cos\theta-1}{\cos\theta+1}\right)=1+ O\big(|r|^{-1}\big),
\end{equation}
and similarly for the other hypergeometric functions. Concentrating on
the integrand of $R_1^{\TX}(\theta)$, in
(\ref{eq:onelevelintegral}), and neglecting factors not depending
on $r$, we find, using~\reff{Fasymptotics}, that the growth-dependence of the integrand on $r$ is
\begin{multline}
2^{2Nr} (1-\cos\theta)^r\prod_{j=0}^{N-1}
\frac{\Gamma(r+1/2+j)}{\Gamma(r+N+j)}\; \frac{\exp(-r\XX)}{r}
\frac{2^{-r}}{2N+r-1} \frac{\Gamma(N+r)}{\Gamma(N+r-1/2)}\times \\
\times
\frac{\Gamma(N+r+1/2)\Gamma(N+r-1/2)}{\Gamma(r+1/2)\Gamma(r-1/2)}
\Big[1+O\big(|r|^{-1}\big)\Big].
\end{multline}
Upon simplification, and replacing the products of Gamma functions with Barnes double
Gamma functions, the growth-dependence on large $|r|$ is
\begin{multline}
 \frac{2^{2Nr-r}}{r(2N+r-1)} \exp(r(\log(1-\cos \theta)-\XX))
\frac{G(r+1/2+N)}{G(r+3/2)} \frac{G(r+N+1)}{G(r+2N)}\times \\
 \times \frac{\Gamma(N+r+1/2)}{\Gamma(r-1/2)}\Big[1+O\big(|r|^{-1}\big)\Big].
\end{multline}

Recall the asymptotic formula for $G(z)$
for large $|z|$:
\begin{equation}
\label{eq:Gasymp}
 \log G(z+1)  \sim z^2(\tfrac{1}{2}\log
  z-\tfrac{3}{4})+\tfrac{1}{2}z\log (2\pi)-\tfrac{1}{12}\log z +\zeta'(-1)+O(z^{-1}).
\end{equation}
Applying~\reff{eq:Gasymp} and Stirling's formula we conclude
that the integrand in (\ref{eq:onelevelintegral}) grows like
\begin{equation} \label{final_growdependence}
\frac{1}{|r|^2}\exp(r(2N\log 2-\log 2 +\log(1-\cos\theta)-\XX))
\Big[1+O\big(|r|^{-1}\big)\Big]
\end{equation}
for large $|r|$.

From~\reff{final_growdependence} we deduce that if $(2N-1)\log 2+\log(1-\cos\theta)-\XX>0$
then we can close the contour in the left half of the complex plane,
thus enclosing the poles at zero and the negative half integers.
If, on the other hand, $(2N-1)\log 2+\log(1-\cos\theta)-\XX <0$
then we must close the contour in the right half-plane, implying that
the integral (\ref{eq:onelevelintegral}) is zero.

Finally, we compute the residues when moving the contour to the left.
There is a simple pole at $r=0$ with residue
\begin{multline}
 2^{N^2-N} \prod_{j=0}^{N-1}
\frac{\Gamma(2+j)\Gamma(1/2+j)\Gamma(1/2+j)}{\Gamma(N+j)}
\frac{2}{2N-1} \frac{\Gamma(N+1)\Gamma(N)}
{\Gamma(N-1/2)\Gamma(N-1/2)}\times \\
 \times P(N,0,\theta).
\end{multline}
Note that when $r=0$ the expression for $P(N,r,\theta)$ reduces to that corresponding to SO$(2N)$.

There is also a simple pole at $r=-1/2$ with residue
\begin{multline}
 -2\exp(\XX/2) 2^{N^2-2N+3/2} \frac{(1-\cos\theta)^{-1/2}} {2N-3/2}
\frac{\Gamma(N+1)\Gamma(1/2)}{\Gamma(N-1)\Gamma(N-1/2)}\times\\
 \times\prod_{j=1}^{N-1} \frac{\Gamma(2+j)\Gamma(1/2+j)\Gamma(j)}
{\Gamma(N+j-1/2)}P(N,-1/2,\theta).
\end{multline}
Similarly to~\reff{CSOasresidues}, where we determined the
normalization constant $\CX$, the contribution from the residues
of the higher order poles at $\alpha_k = \frac{-1-2k}{2}$ with
$k=1,2,3,\ldots$ is
\begin{equation} \label{eq:residuesum}
\sum_{k \geq 1}^\infty b_k \exp((k+1/2)\XX)
\end{equation}
where the coefficients $b_k$ come from the residues.
\end{proof}
To show convergence we may write the series~\reff{eq:residuesum} as
\begin{align}
\sum_{1 \leq k \leq K} b_k \exp((k+1/2)\XX) + \int_{c_K-i\infty}^{c_K+i\infty} \exp(-r\XX)f(r) dr
\end{align}
where $\exp(-r\XX)f(r)$ is the integrand in~\reff{eq:onelevelintegral} and $\alpha_{K+1} < c_K < \alpha_K$. By~\reff{final_growdependence},
$f(r)\ll \exp(-rX)$ for large $|r|$ and $\XX < 0$. Thus,
\begin{equation}
  \begin{split}
    \int_{c_K-i\infty}^{c_K+i\infty}\exp(-r\XX) f(r) dr & = \exp(-c_K \XX) i \int_{-\infty}^{\infty} \exp(-it\XX) f(c_K+it) dt \\
    & = O(\exp(-c_K \XX)).
  \end{split}
\end{equation}
Hence, the contributions of residues of higher poles decrease
exponentially fast for $\XX < 0$ and using only the $K$ right-most
poles gives an error that is $O(\exp(-c_K \XX))$. The computation
of residues of higher order poles can be carried out by computer
algebra.  We find that for reasonable $\XX$ we have good
convergence using about 10 poles.

\begin{remark} \rm The hard gap reflected in Theorem
 ~\ref{corollaryone} may be understood as follows. Recalling that
  $\Lambda_A(1,N)=2^N\prod_{n=1}^N(1-\cos\theta_n)$ it follows that
  $\log \Lambda_A(1,N) = (2N-1)\log2 +\log(1-\cos\theta)$ for the
  matrix $A\in\text{SO}(2N)$ having $2N-2$ eigenvalues at $-1$ (i.~e.,
  $\theta_2,\dots,\theta_N=\pi$) and a symmetric pair of eigenvalues
  at $e^{\pm i\theta}$ (i.~e., $\theta_1=\theta$).  Plainly, such
  matrix~$A$ maximizes $\log\Lambda_A(1,N)$ subject to the condition
  that a pair of eigenvalues lie at $e^{\pm i\theta}$. Therefore, for
  any fixed $\XX$, the condition $\log\Lambda_A(1,N)>\XX$ implies
  $\XX<(2N-1)\log2 + \log(1-\cos\theta)$: if the latter inequality
  does not hold, there are \emph{no} matrices~$A$ with
  $\log\Lambda_A(1,N)>\XX$. In particular, for fixed $N$ and $\XX$
  there is a lower bound for all eigenphases $\theta$ of any matrix
  $A$ with $\log\Lambda_A(1,N)>\XX$, namely $\theta >
  \theta_{\text{inf}}:= \cos^{-1}(1-2^{-(2N-1)}e^{\XX})$: the excised
  one-level density is identically zero in the excluded spectral
  interval (``hard gap'') $\theta\le\theta_{\text{inf}}$. Note that
  the hard gap shrinks to zero exponentially fast as $N\to\infty$ and
  also as $\XX \rightarrow -\infty$.
\end{remark}

\subsection{Formula vs. data}
\begin{center}
\begin{figure}[h]
\includegraphics[scale=0.5,angle=-90]{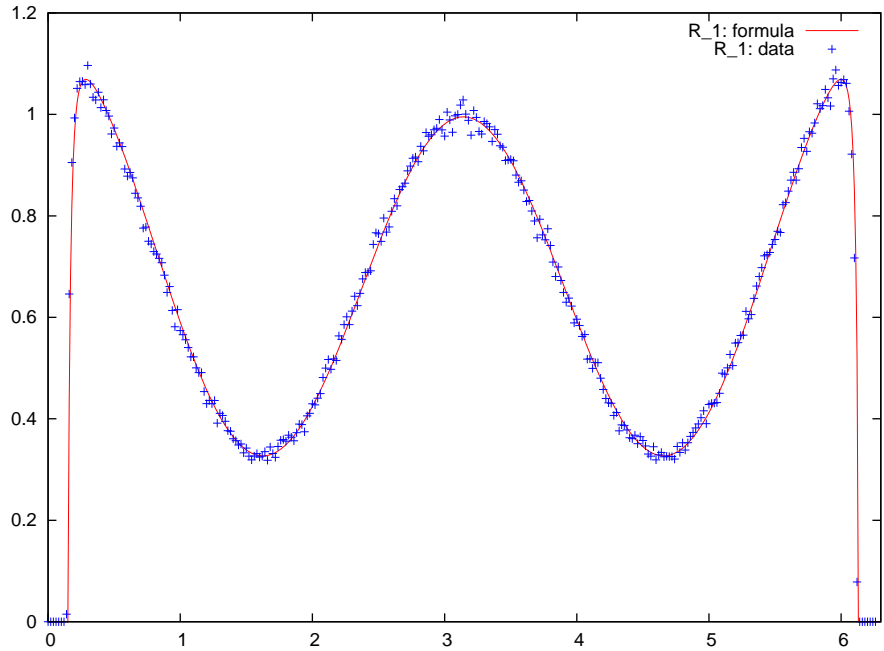}
\caption{\label{fig:consistency} One-level density of excised
SO($2N), N = 2$
  with cut-off $|\Lambda_A(1,N)| \geq 0.1$. The \emph{red curve} uses our
  formul\ae\ from Theorems~\ref{theoremone} and~\ref{corollaryone}. The
  \emph{blue crosses} give the empirical one-level density of 200,000
  numerically generated matrices.}
\end{figure}
\end{center}

For illustration and as a consistency test, in figure~\ref{fig:consistency} we compare our formula for the one-level density
with data obtained by generating random matrices from the excised
ensemble SO($2N)$.

\newcommand{\SOeven}{{\rm SO}(2N)}

\bibliography{finiteconductor}{}
\bibliographystyle{plain}

\ \\

\end{document}